\documentclass[12pt]{article}
\usepackage[russian]{babel}
\usepackage{amsthm,amsfonts,amssymb,amscd}

\begin{document}

\begin{center}
\Large \bf Birationally rigid \\
complete intersections of quadrics and cubics
\end{center}
\vspace{1cm}

\centerline{A.V.Pukhlikov}

\parshape=1
3cm 10cm \noindent {\small \quad \quad \quad
\quad\quad\quad\quad\quad\quad\quad {\bf }\newline We prove
birational superrigidity of generic Fano complete intersections
$V$ of type $2^{k_1}\cdot 3^{k_2}$ in the projective space
${\mathbb P}^{2k_1+3k_2}$, under the condition that $k_2\geq 2$
and $k_1+2k_2=\mathop{\rm dim} V\geq 12$, and of a few families of
Fano complete intersections of dimension 10 and 11.

Bibliography: 24 titles.} \vspace{1cm}

{\large \bf Introduction}\vspace{1cm}

{\bf 0.1. Formulation of the main result.} Let $k_1,k_2\in{\mathbb
Z}_+$ be a pair of nonnegative integers, $k_1+2k_2\geq 4$. A
smooth $(k_1+2k_2)$-dimensional complete intersection
$$
V=F_1\cap\dots\cap F_{k_1+k_2}\subset{\mathbb P}={\mathbb P}^{2k_1+3k_2}
$$
of codimension $k=k_1+k_2$ in the projective space ${\mathbb P}$
over the field of complex numbers ${\mathbb C}$, where
$F_i\subset{\mathbb P}$ are hypersurfaces of degree 2 for
$i=1,\dots,k_1$ and of degree 3 for $i=k_1+1,\dots,k=k_1+k_2$ will
be called  a Fano complete intersection of the type $2^{k_1}\cdot
3^{k_2}$. Set $M=k_1+2k_2=\mathop{\rm dim}V$. Obviously, $V$ is a
Fano variety of index one, $\mathop{\rm Pic}V={\mathbb Z}H$ and
$K_V=-H$, where $H$ is the class of a hyperplane
section.\vspace{0.1cm}

The main result of this paper is the following\vspace{0.1cm}

{\bf Theorem 1.} {\it Assume that $k_2\geq 2$ and either
$M=k_1+2k_2\geq 12$, or the pair $(k_1,k_2)$ is one of the
following five pairs:
$$
(5,3),(3,4), (1,5), (2,4), (0,5).
$$
Then a generic (in the sense of the Zariski topology on the space
of coefficients of the polynomials, defining the hypersurfaces
$F_1,\dots,F_{k_1+k_2}$) Fano complete intersection $V$ of type
$2^{k_1}\cdot 3^{k_2}$ is birationally superrigid. In
particular,\vspace{0.1cm}

{\rm (i)} there exists no rational dominant map $\gamma\colon
V\dashrightarrow S$ onto a variety $S$ of positive dimension, the
generic fibre of which $\gamma^{-1}(s)$ is of negative Kodaira
dimension,\vspace{0.1cm}

{\rm (ii)} any birational map $\chi\colon V\dashrightarrow V'$
onto a Fano variety $V'$ with ${\mathbb Q}$-factorial terminal
singularities and the Picard number $\mathop{\rm rk}\mathop{\rm
Pic}V'=1$ is a biregular isomorphism,\vspace{0.1cm}

{\rm (iii)} the group of birational automorphisms $\mathop{\rm
Bir}V$ coincides with the group of biregular (projective)
automorphisms $\mathop{\rm Aut}V$ and for that reason is
trivial.}\vspace{0.1cm}

The properties (i)-(iii) follow immediately from the birational
superrigidity, see~\cite{Pukh07b}. In its turn, the birational
superrigidity (understood in the sense of equality of the virtual
and actual thresholds of canonical adjunction,
$$
c_{\rm virt}(\Sigma)=c(\Sigma)=n,
$$
for any mobile linear system $\Sigma\subset|-nK_V|$, see
\cite{Pukh07b} for the definitions) follows immediately from the
canonicity of every pair $(V,\frac{1}{n}\Sigma)$, where
$\Sigma\subset|-nK_V|$ is a mobile system, or, equivalently, from
the fact that the linear system $\Sigma$ has no maximal
singularities.\vspace{0.1cm}

Up to this day, birational superrigidity was shown for Fano
complete intersections of index one $F_1\cap\dots\cap F_k$ under
the condition that for at least one degree the inequality
$\mathop{\rm deg}F_i\geq 4$ holds \cite{Pukh01,Pukh11}. Theorem 1
states the superrigidity for complete intersections $V$ of
quadrics and cubics under the assumption that $\mathop{\rm
dim}V\geq 12$ and there are at least two cubics and also for
certain five families of complete intersections of dimension 10
and 11. Thus the problem of birational superrigidity (of a generic
variety) remains open for several families of Fano complete
intersections of dimension $\leq 11$ and for the following three
infinite series:
$$
2\cdot 2\cdot\dots\cdot 2,\quad 2\cdot 2\cdot\dots\cdot 2\cdot 3
\quad \mbox{and}\quad 2\cdot 2\cdot\dots\cdot 2\cdot 4
$$
in the arbitrary dimension. Those varieties require further
improvement of the techniques of the proof.\vspace{0.1cm}

In the next subsection we give a precise meaning to the condition
of general position of the variety $V$ in its family. The main
difficulty in the proof of Theorem 1 is to show that this
condition does realize, that is, there is a non-empty Zarisky open
subset in the space of parameters, corresponding to the varieties
$V$, satisfying that condition. The very proof of birational
superrigidity take a few pages (\S1), almost all paper is devoted
to the proof of the conditions of general position.\vspace{0.1cm}

In Sec. 0.3 we give a more detailed description of the structure
of the paper, of the main idea of the proof of Theorem 1 and of
the position of that theorem among the other results on birational
rigidity.\vspace{0.3cm}

{\bf 0.2. Conditions of general position.} As usual, the
conditions of general position (or the {\it regularity
conditions}) are of local nature and should be satisfied at every
point of the variety $V$. Let $o\in V\subset{\mathbb P}$ be some
point, ${\mathbb C}^{M+k}\subset{\mathbb P}$ a standard affine set
with the coordinates $(z_1,\dots,z_{M+k})$, where $o=(0,\dots,0)$.
Let
$$
f_i=q_{i,1}+q_{i,2}
$$
for $i=1,\dots,k_1$, be the equations of the quadrics
$F_1,\dots,F_{k_1}$ and
$$
f_i=q_{i,1}+q_{i,2}+q_{i,3}
$$
for $i=k_1+1,\dots,k_1+k_2$, the equations of the cubics
$F_{k_1+1},\dots,F_k$, decomposed into homogeneous in $z_*$
components. The tangent space $T_oV\subset{\mathbb C}^{M+k}$ is
given be the system of equations
$$
q_{1,1}=\dots=q_{k,1}=0.
$$
We denote its projectivization ${\mathbb P}(T_oV)\cong{\mathbb
P}^{M-1}$ by the symbol ${\mathbb T}$. Set
$\bar{q}_{i,j}=q_{i,j}|_{\mathbb T}$. Let us introduce the
following sets of pairs of indices:
$$
J_1=\{(i,2)\,|\,1\leq i\leq k\},\quad J_2=\{(i,3)\,|\,k_1+1\leq
i\leq k_1+k_2\}
$$
and $J=J_1\cup J_2$. The first (traditional) regularity condition
is formulated in the following way:\vspace{0.1cm}

(R1) The system of homogeneous equations
\begin{equation}\label{9april12.5}
\{\bar{q}_{i,2}=0\,|\,1\leq i\leq k\}
\end{equation}
defines in ${\mathbb T}$ a closed subset of pure codimension $k$,
whereas the system of equations
\begin{equation}\label{9april12.6}
\{\bar{q}_{i,j}=0\,|\,1(i,j)\in J\}
\end{equation}
defines in ${\mathbb T}$ either the empty set or a finite set of
points, linearly independent in ${\mathbb T}$.\vspace{0.1cm}

The last condition (on the linear independence) is a new addition
to the usual regularity condition \cite{Pukh01,Pukh11}, which
requires only that the set (\ref{9april12.6}) should be finite. It
is clear that if the condition (R1) is satisfied at the point $o$,
then on $V$ lie at most $M$ lines, passing through that point. The
second regularity condition was introduced in
\cite{Pukh11}.\vspace{0.1cm}

(R2) None of the irreducible components of the closed algebraic
set (\ref{9april12.5}) is contained in a linear subspace of
codimension two in ${\mathbb T}$.\vspace{0.1cm}

In \cite{Pukh11} the condition (R2) was called {\it correctness in
the quadratic terms}.\vspace{0.1cm}

To formulate the third regularity condition, we need some
additional constructions. Assume that the set (\ref{9april12.6})
is non-empty and denote it by the symbol $\Delta$. Let
$$
\Sigma^{\mathbb T}_1=<\bar{q}_{i,2}\,|\,1\leq i\leq k>
$$
be the linear system of quadrics on ${\mathbb T}$, generated by
$\bar{q}_{i,2}$ and $\Sigma^{\mathbb T}_2$ the linear system of
cubics, consisting of all linear combinations of the form
$$
\sum^k_{i=1}
\bar{s}_i(z)\bar{q}_{i,2}+\sum^{k_1+k_2}_{i=k_1+1}\lambda_i\bar{q}_{i,3},
$$
where $\lambda_i$ are constants and $\bar{s}_i(z)=s_i(z)|_{\mathbb
T}$ are linear polynomials on ${\mathbb T}$. Let
$$
{\cal P}=\underbrace{\Sigma^{\mathbb T}_1\times\dots
\times\Sigma^{\mathbb T}_1}_{k-2}
\times\underbrace{\Sigma^{\mathbb T}_2\times\dots
\times\Sigma^{\mathbb T}_2}_{k_2-2}
$$
be the space of all tuples of polynomials $(g_1,\dots,g_{M-4})$,
where the first $k-2$ polynomials are taken from the space
$\Sigma^{\mathbb T}_1$, and the last $k_2-2$ from the space
$\Sigma^{\mathbb T}_2$. Now let us formulate the third and last
regularity condition:\vspace{0.1cm}

(R3) for any irreducible subvariety $R\subset{\mathbb T}$ of
degree $d_R\geq 1$ and codimension 3 there exists a non-empty
Zarisky open subset $U_R\subset{\cal P}$, such that for any tuple
$$
(g_1,\dots,g_{M-4})\in U_R
$$
the following inequality holds:
\begin{equation}\label{9april12.7}
\sum_{p\in \Delta}\mathop{\rm dim}{\cal O}_{p,R}/(g_1,\dots,g_{M-4})
\leq 2^{k-4}3^{k_2-1}d_R.
\end{equation}
Note that by the condition (R1) for a generic tuple
$(g_1,\dots,g_{M-4})\in{\cal P}$ the scheme
$$
(R\circ\{g_1=0\}\circ\dots\circ\{g_{M-4}=0\})
$$
is zero-dimensional and for that reason is of degree
$2^{k-2}3^{k_2-2}d_R$. The condition (R3) means that at most $3/4$
of that full degree is concentrated at the points of the set
$\Delta$. This is a very strong condition, since it should be
satisfied for an {\it arbitrary} subvariety $R$ of the projective
space ${\mathbb T}$.\vspace{0.1cm}

Let
$$
{\cal F}_{sm}\subset \prod^{k_1}_{i=1}H^0({\mathbb P},{\cal
O}_{\mathbb P}(2))\times\prod^k_{i=k_1+1}H^0({\mathbb P},{\cal
O}_{\mathbb P}(3))
$$
be the space of tuples $(f_1,\dots,f_k)$, defining a smooth
complete intersection of the type $2^{k_1}3^{k_2}$ in ${\mathbb
P}$.\vspace{0.1cm}

{\bf Theorem 2.} {\it There exists a non-empty Zarisky open subset
${\cal F}_{\rm reg}\subset {\cal F}_{\rm sm}$ of tuples
$(f_1,\dots,f_k)$, for which the corresponding complete
intersection $V(f_1,\dots,f_k)=\{f_1=\dots=f_k=0\}\subset{\mathbb
P}$ satisfies the regularity conditions (R1)-(R3) at every point
$o\in V$.}\vspace{0.1cm}

Now we can make the main result of the present paper more
precise.\vspace{0.1cm}

{\bf Theorem 3.} {\it Let $V$ be a Fano complete intersection of
the type $2^{k_1}3^{k_2}$ in ${\mathbb P}$, satisfying the
regularity conditions (R1)-(R3) at every point $o\in V$, where the
numbers $k_1,k_2$ satisfy the assumptions of Theorem 1. Then the
variety $V$ is birationally superrigid. In particular, the claims
(i)-(iii) of Theorem 1 hold.}\vspace{0.1cm}

Obviously, Theorem 1 follows from Theorems 2 and 3, which are
independent of each other and will be shown
separately.\vspace{0.3cm}

{\bf 0.3. The structure of the paper and historical remarks.}
Theorem 1 is proved in \S1. Its proof takes a few pages and makes
use of the technique of hypertangent linear system and the
$8n^2$-inequality. In two words, the idea of the proof can be
explained as follows. In the preceding papers (see [1-5] and the
bibliography in \cite{Pukh07b}) to show birational superrigidity
of Fano complete intersections, one constructed, starting from the
self-intersection $Z$ of a movable linear system $\Sigma$ on $V$
with a maximal singularity, an effective1-cycle $Y\subset V$,
satisfying the inequality
\begin{equation}\label{9april12.8}
\mathop{\rm mult}\nolimits_oY>\mathop{\rm deg}Y,
\end{equation}
which is, of course, impossible, and brings the assumption that
such a system $\Sigma$ exists to a contradiction, thus proving
birational superrigidity of the variety $V$. For complete
intersections of quadrics and cubics this construction does not
work in the straightforward way: the technique of hypertangent
systems makes it possible to construct a curve $Y\subset V$ that
has a high multiplicity $\mathop{\rm mult}_oY$, which, however,
does not exceed its degree $\mathop{\rm deg}Y$, so that there is
no contradiction.\vspace{0.1cm}

An easy analysis shows that the irreducible curve $C\subset V$,
{\it which is not a line} (that is, a curve of degree $\mathop{\rm
deg}C\geq 2$), satisfies an essentially stronger inequality than
the inequality, opposite to (\ref{9april12.8}): the estimate
\begin{equation}\label{9april12.9}
\mathop{\rm mult}\nolimits_o C\leq\frac23\mathop{\rm deg}C
\end{equation}
holds. Breaking the effective 1-cycle $Y$ into two sub-cycles:
$$
Y=Y_1+Y_{\geq 2},
$$
where $Y_1$ is concentrated on the lines, passing through the
point $o$, and $Y_{\geq 2}$ is concentrated on the curves of
degree $\geq 2$, we see that the strategy of the proof, described
above, would be successful if the sub-cycle $Y_1$ {\it is not too
large}, that is, if the ratio $\mathop{\rm deg}Y_1/\mathop{\rm
deg}Y$ is small enough. Using the more precise estimate
(\ref{9april12.9}) for the sub-cycle $Y_{\geq 2}$, one can get a
contradiction and prove birational rigidity. This is possible,
provided that the ``outcome'' of the technique of hypertangent
divisors, applied to the self-intersection $Z$ of the linear
system $\Sigma$ (or to a certain irreducible component of that
self-intersection), is an effective 1-cycle $Y$, containing ``not
too many lines''. Proof of Theorem 3, realizing the strategy,
described above, is given in \S1.\vspace{0.1cm}

The main part of the paper, \S\S2-5, is about estimating the
number of lines (taking into account the multiplicities), emerging
as the result of intersecting an arbitrary irreducible subvariety
with hypertangent divisors. In \S\S2-4 we consider the local
problem: to estimate the multiplicity of a generic tuple of
polynomials from a given linear system with an arbitrary effective
cycle at a fixed point. In \S5 the problem is globalized: we
estimate the sum of local intersection multiplicities over the
base points of the linear system, and on that basis prove Theorem
2.\vspace{0.1cm}

The theory of birational rigidity of Fano complete intersections
has a long history, starting from the work of Fano himself
\cite{Fano3}, where he studied complete intersections $V_{2\cdot
3}\subset{\mathbb P}^5$ of a quadric and a cubic and complete
intersections of three quadrics $V_{2\cdot 2\cdot
2}\subset{\mathbb P}^6$. Fano gave a description (as it turned out
later, non-complete \cite{I80}) of generators of the group of
birational self-maps $\mathop{\rm Bir}V_{2\cdot 3}$ and formulated
the theorem on non-rationality for the both classes of varieties.
As it was discovered later \cite{IM}, his arguments contained
numerous mistakes and gaps and can not be considered even as a
first approximation to the rigorous proof. At the same time, in
his works Fano outlined many important ideas and constructions
which later were proved essential in birational geometry (such as
the ``double projection'' or the Noether-Fano
inequality).\vspace{0.1cm}

In 1970 in their pioneer paper \cite{IM} V.A.Iskovskikh and
Yu.I.Manin gave the first ever rigorous proof of birational
superrigidity (in the modern terminology) for one class of
three-dimensional Fano varieties, the smooth three-dimensional
quartics $V_4\subset{\mathbb P}^4$. With a minimal modification
(which deal with the easier part of the proof, that on exclusion
or untwisting maximal curves lying on the variety itself) the
proof of V.A.Iskovskikh and Yu.I.Manin worked for the Fano double
spaces $V_2\stackrel{2:1}{\to}{\mathbb P}^3$ of index 1 and the
double quadrics $V_4\stackrel{2:1}{\to}Q_2\subset{\mathbb P}^4$ of
index 1, see \cite{I80}. However, the proof of birational rigidity
and description of the group of birational self-maps of the
variety $V_{2\cdot 3}$ (announced in \cite{I77b}) turned out to be
a much harder problem. The {\it test class method}, developed in
\cite{IM}, worked successfully for varieties of degree at most 4
only, and already the degree $\mathop{\rm deg}V_{2\cdot 3}=6$ made
a serious obstruction. The arguments of \cite{I80}, used for the
exclusion of the infinitely near maximal singularity \cite[Sec.
4.8, Step 3]{I80}, were erroneous, and this, most essential step
of the proof remained an open problem until 1987 \cite{Pukh89b}. A
complete proof of the theorem on birational rigidity of a generic
complete intersection of a quadric and a cubic of dimension three
see in \cite{IP}. The problem of description of the structures of
a rationally connected fibre space and the group of birational
self-maps of the variety $V_{2\cdot 2\cdot 2}$ is still open.
(Note that birational geometry of Fano complete intersections was
also studied by means of the transcendental method; there are not
many such papers, see \cite{Tyu75b,Tyu79}.)\vspace{0.1cm}

The work on extending the test class method into arbitrary
dimension was started in \cite{Pukh87,Pukh89a}. The ideas
developed in \cite{Pukh87} are sufficient for the proof of
birational superrigidity of generic complete intersections of a
quadric and a quartic $V_{2\cdot 4}\subset{\mathbb P}^6$ of
dimension 4 and index 1. This result was announced in 1985, but
its complete proof was not published. (Later in \cite{Ch03}
birational superrigidity of complete intersections of a quadric
and a quartic, not containing planes, was shown via the
$8n^2$-inequality, however, the proof of that inequality that was
known at the time turned out to be erroneous and a complete proof
was obtained later, see the history of this problem in
\cite{Pukh10}.) \vspace{0.1cm}

Replacing the test class method by the method of counting
multiplicities, introdu\-cing the techniques of hypertangent
divisors and some other ideas \cite{Pukh98b} made it possible to
prove birational superrigidity of Fano hypersurfaces of index 1
and later of generic complete intersections $V_{d_1\cdot\dots\cdot
d_k}\subset{\mathbb P}^{M+k}$ of index 1 for $M\geq 2k+1$,
$\mathop{\rm dim}V\geq 4$ \cite{Pukh01}. The next step in the
development of this area was systematic use of the connectedness
principle of Shokurov and Koll\' ar \cite{Kol93} and the proof of
the $8n^2$-inequality, which formed a basis for proving birational
superrigidity for a wider class of complete intersections with
$M\geq k+3$, $d_k\geq 4$ \cite{Pukh11}. The latter paper is the
immediate predecessor of the present one: here we remove the
restriction $d_k\geq 4$. As we noted above, now the problem of
birational (super)rigidity remains open for very few classes of
Fano complete intersections of index 1. The conjecture on
bitational superrigidity of arbitrary smooth Fano complete
intersections of index 1 was formulated by the author many times
(see, for instance, \cite{Pukh04c}).\vspace{0.1cm}

The method of maximal singularities makes it also possible to
prove birational superrigidity of the varieties fibred into Fano
complete intersections over the projective line ${\mathbb P}^1$,
see \cite{Pukh06b,Pukh09a}, but this is a different
topic.\vspace{0.1cm}

The work on the problem, a solution of which makes the contents of
the present paper, was completed in spring 2012. As far as the
author knew, the problem of estimating the multiplicity for a
generic tuple of polynomials in a subvariety of a given
codimension has never been considered before. However, in June
2012 A.G.Khovanskii informed the author that a different, but
close problem, that of estimating the multiplicity of an isolated
solution of a system of $n$ equations
$$
\phi_1=\dots=\phi_n=0
$$
in $n$ complex variables for a tuple $(\phi_1,\dots,\phi_n)$ in a
ring of Noetherian functions $K\supset {\mathbb
C}[z_1,\dots,z_n]$, finitely generated over ${\mathbb C}[z_*]$ and
closed with respect to differentiations, was studied in
\cite{GabKh}. In the latter paper the estimate was obtained in
terms of invariants of the ring $K$ using the methods that were
essentially different from the algebro-geometric techniques used
in \S\S 2-4: by reduction to the one-dimensional problem of
restricting a polynomial onto a trajectory of a polynomial vector
field, via estimating the complexity of an integral manifold of an
analytic vector-function and using deformations. In our paper we
consider a less general (and in fact somewhat different) problem,
but a considerably stronger and (which is especially important for
applications to birational geometry) effective estimate of the
multiplicity, which is not implied by the estimates of
\cite{GabKh}.\vspace{0.1cm}

The method of solution, developed in this paper (that is, the
method of the proof of Theorem 2), is new. A simplified version of
this method (for a system of $N$ polynomial equations in $N$
variables, without involving an effective cycle $R$ of the given
degree) was published in \cite{Pukh12}. \vspace{0.1cm}

The very fact that the problem of estimating multiplicity of an
(isolated) solution of a system of $n$ equations in $n$ variables
emerges in different contexts, requires different types of
techniques and has various applications, is remarkable; the author
is grateful to A.G.Khovanskii for pointing out the paper
\cite{GabKh}.\vspace{0.5cm}

%%%%%%%%%%%%%%%%%%%%%%%%%%%%%%%%%%%%%%%%%%%%%%%%%%%%%%%%%%%%%%%%%%
%%%%%%%%%%%%%%%%%%%%%%%%%%%%%%%%%%%%%%%%%%%%%%%%%%%%%%%%%%%%%%%%%%
%%%%%%%%%%%%%%%%%%%%%%%%%%%%%%%%%%%%%%%%%%%%%%%%%%%%%%%%%%%%%%%%%%
%%%%%%%%%%%%%%%%%%%%%%%%%%%%%%%%%%%%%%%%%%%%%%%%%%%%%%%%%%%%%%%%%%
%%%%%%%%%%%%%%%%%%%%%%%%%%%%%%%%%%% SECTION 1

{\large\bf \S 1. Proof of birational superrigidity}\vspace{0.3cm}

In this section we prove Theorem 3.\vspace{0.3cm}

{\bf 1.1. Start of the proof. The $8n^2$-inequality.} The proof of
Theorem 3 starts in the standard way. Let us fix a complete
intersection $V\subset{\mathbb P}$, satisfying the assumptions of
the theorem, and a mobile linear system $\Sigma\subset |nH|$
(where, let us remind the reader, $H$ is the class of a hyperplane
section of the variety $V$) with a {\it maximal singularity}
\cite{Pukh07b}. The latter means that for some birational morphism
$\varphi\colon V^{\sharp}\to V$, where $V^{\sharp}$ is a
non-singular projective variety, and an exceptional divisor
$E^{\sharp}\subset V^{\sharp}$ (the maximal singularity) the {\it
Noether-Fano inequality} holds
$$
\mathop{\rm ord}\nolimits_{E^{\sharp}}\varphi^{*}\Sigma>na(E^{\sharp},V).
$$
The irreducible subvariety $B=\varphi(E^{\sharp})\subset V$ is
called the {\it centre} of the maximal singularity $E^{\sharp}$.
The inequality $\mathop{\rm mult}_B\Sigma>n$ holds. There are
three options for the codimension of the subvariety
$B$:\vspace{0.1cm}

(i) $\mathop{\rm codim}B=2$,\vspace{0.1cm}

(ii) $\mathop{\rm codim}B=3$,\vspace{0.1cm}

(iii) $\mathop{\rm codim}B\geq 4$.\vspace{0.1cm}

As it was shown in, for example, \cite{Pukh01} (that argument can
be found in many papers on birational rigidity), the first option
does not realize, since for the {\it numerical} Chow group we have
$A^2V={\mathbb Z}H^2$. The second option is excluded, for
instance, in \cite{Pukh11}. Thus we will assume that $\mathop{\rm
codim}B\geq 4$.\vspace{0.1cm}

At first our arguments repeat the proof of Theorem 3 in
\cite[Sec.~3]{Pukh11} almost word for word. We assume that the
codimension of the subvariety $B$ is minimal among all centres of
maximal singularities of the system $\Sigma$; in particular, $B$
is not strictly contained in the centre of another maximal
singularity, if there are any. Take a point of general position
$o\in B$. Let $\lambda\colon V^+\to V$ be its blow up,
$E=\lambda^{-1}(o)\cong{\mathbb P}^{M-1}$ the exceptional divisor,
$\Sigma^+$ the strict transform of the mobile system $\Sigma$ on
$V^+$, $Z=(D_1\circ D_2)$ and $Z^+$ are the self-intersection of
the system $\Sigma$ and its strict transform on $V^+$,
respectively.\vspace{0.1cm}

{\bf Proposition 1.1 (the $8n^2$-inequality).} {\it There exist a
linear subspace $P\subset E$ of codimension two, satisfying the
inequality
$$
\mathop{\rm mult}\nolimits_oZ+\mathop{\rm mult}\nolimits_PZ^+>8n^2.
$$
If $\mathop{\rm mult}_o Z\leq 8n^2$, then the linear subspace $P$
is uniquely determined by the system} $\Sigma$.\vspace{0.1cm}

{\bf Proof} was given in \cite[Sec.~4.1]{Pukh10}.\vspace{0.1cm}

Let $\Lambda\supset P$ be a generic hyperplane in $E\cong{\mathbb
P}^{M-1}$, containing the subspace $P$, $L\in |H|$ a generic
hyperplane section, containing the point $o$ and such that
$L^+\cap E=\Lambda$, where $L^+$ is the strict transform of the
divisor $L$ on $V^+$. By genericity of the choice of $\Lambda$ and
$L$ none of the irreducible components of the cycle $Z$ is
contained in $L$, so that the scheme-theoretic intersection
$Z_L=(Z\circ L)$ is well defined. The effective cycle $Z_L$ of
codimension three satisfies the inequality
$$
\mathop{\rm mult}\nolimits_oZ_L\geq\mathop{\rm mult}\nolimits_oZ+
\mathop{\rm mult}\nolimits_PZ^+>8n^2
$$
and the equality $\mathop{\rm deg}Z_L=\mathop{\rm deg}Z=dn^2$,
where $d=2^{k_1}3^{k_2}$ is the degree of the complete
intersection $V$. Therefore, there exists an irreducibel
subvariety $Q\subset V$ of codimension three (an irreducible
component of the cycle $Z_L$), satisfying the estimate
\begin{equation}\label{9april12.1}
\frac{\mathop{\rm mult}\nolimits_oQ}{\mathop{\rm deg}Q}>\frac{8}{d}.
\end{equation}

Let $Q_E=(Q^+\circ E)=\sum\limits_{i\in I}m_iR_i$ be the
projectivized tangent cone to the subvariety $Q$ at the point $o$
(here $Q^+\subset V^+$ is the strict transform of the subvariety
$Q$, the subvarieties $R_i\subset E$ are irreducible components of
the effective cycle $Q_E$, $m_i\geq 1$; here and everywhere in
this paper constructions of elementary intersection theory are
understood in the sense of \cite{Ful}). Now let us apply to the
subvariety $Q\subset V$ the technique of hypertangent divisors in
the way in which it was done in \cite{Pukh11}, but with a small
modification. Namely, at each step when intersecting with a
(hyper)tangent divisor, we will ignore only the emerging
irreducible components that do not contain the point
$o$.\vspace{0.3cm}

%%%%%%%%%%%%%%%%%%%%%%%%%%%%%%%%%%%%%%%%%%%%%%%%%%%%%%%%%%%%%%%%%%%%%%%%%%%%%%%%%%%%
%%%%%%%%%%%%%%%%%%%%%%%%%%%%%%%%%%%%%%%%%%%%%%%%%%%%%%%%%%%%%%%%%%%%%%%%%%%%%%%%%%%%
%%%%%%%%%%%%%%%%%%%%%%%%%%%%%%%%%%% subsection 1.2

{\bf 1.2. The technique of hypertangent divisors.} Let
$$
\Sigma_1=\{\lambda_1q_{1,1}|_V+\dots+\lambda_kq_{k,1}|_V\}
$$
be the $k$-dimensional space of equations of tangent hyperplanes
at the point $o$. By the condition (R1) for any non-zero tuple
$(\lambda_1,\dots,\lambda_k)$ the corresponding tangent divisor
$$
T(\lambda_*)=\{\lambda_1q_{1,1}|_V+\dots+\lambda_kq_{k,1}|_V=0\}
$$
satisfies the equality
$$
T^E(\lambda_*)=(T^+(\lambda_*)\circ E)=
\{\lambda_1\bar{q}_{1,2}+\dots+\lambda_k\bar{q}_{k,2}=0\}.
$$
In particular, set $T_i=\{q_{i,1}|_V=0\}$. Since none of the
irreducible components of the effective cycle
$(T^E_1\circ\dots\circ T^E_k)$ of codimension $k$ by the condition
(R2) is contained in the linear subspace of codimension two in
$E$, we get the equality
\begin{equation}\label{9april12.2}
\mathop{\rm codim}\nolimits_E(T^E_1\cap\dots\cap T^E_k\cap\Lambda)=k+1.
\end{equation}
The support of the cycle $Q_E$, constructed above, is contained in
the hyperplane $\Lambda$. By the equality (\ref{9april12.2}), for
$(k-2)$ generic divisors $D_{1,1},\dots,D_{1,k-2}$ in the linear
system $\Sigma_1$ the equality
$$
\mathop{\rm codim}\nolimits_E(Q_E\cap D^E_{1,1}\cap\dots\cap
D^E_{1,k-2})=k+1
$$
holds, where $D^E_{1,i}=(D^+_{1,i}\circ E)$. Therefore, in a
neighborhood of the point $o$ we also have
\begin{equation}\label{9april12.3}
\mathop{\rm codim}\nolimits_o(Q\cap D_{1,1}\cap\dots\cap
D_{1,k-2})=k+1.
\end{equation}
Now let us construct the following sequence of effective cycles
$Q_i$, where $i=0,1.\dots,k-2$:\vspace{0.1cm}

(i) $Q_0=Q$,\vspace{0.1cm}

(ii) $Q_{i+1}$ is obtained by removing from the effective cycle
$(Q_i\circ D_{1,i+1})$ all irreducible components, not containing
the point $o$.\vspace{0.1cm}

This procedure is well defined: since all components of each of
the cycles $Q_i$ contain the point $o$, by the equality
(\ref{9april12.3}) none of the components of the cycle $Q_i$ is
contained in the divisor $D_{1,i+1}$. In particular, $\mathop{\rm
codim}Q_i=i+3$ for $i=0,\dots,k-2$. The degree $\mathop{\rm
deg}Q_i$ is not increasing (in fact, it is decreasing, if some
components are indeed removed in the process of this
construction), so that $\mathop{\rm deg}Q_{k-2}\leq\mathop{\rm
deg}Q$, whereas for the multiplicity at the point $o$ we get the
equality
$$
\mathop{\rm mult}\nolimits_oQ_{k-2}=
2^{k-2}\mathop{\rm mult}\nolimits_oQ.
$$
Now let us consider the hypertangent linear system
$$
\Sigma_2=\{h(z)|_V=\sum^k_{i=1}s_i(z)q_{i,1}|_V+
\sum^{k_1+k_2}_{i=
k_1+1}\lambda_i(q_{i,1}+q_{i,2})|_V\}.
$$
Obviously, $\Sigma^+_2\subset|2H-3E|$: the projectivized tangent
cone of the divisor $h|_V=0$ at the point $o$ is given by the
equation
$$
-\left(\sum^k_{i=1}\bar{s}_i\bar{q}_{i,2}+
\sum^{k_1+k_2}_{i=k_1+1}\lambda_i\bar{q}_{i,3}\right)\in\Sigma^{\mathbb
T}_2.
$$
Let $(D_{2,1},\dots,D_{2,k_2-2})\in\Sigma^{\times(k_2-2)}_2$ be a
generic tuple of $(k_2-2)$ hypertangent divisors. By the condition
(R1) the closed set
$$
Q_E\cap\left(\mathop{\bigcap}\limits^{k-2}_{i=1}D^E_{1,i}\right)\bigcap
\left(\mathop{\bigcap}\limits^{k_2-2}_{i=1}D^E_{2,i}\right)
$$
is zero-dimensional, so that the closed set
$$
Q\cap\left(\bigcap^{k-2}_{i=1}D_{1,i}\right)\cap
\left(\bigcap^{k_2-2}_{i=1}D_{2,i}\right)
$$
is one-dimensional in a neighborhood of the point $o$. Let us
continue to construct the chain of effective cycles $Q_i$,
$i=k-2,\dots,M-4$, where $Q_{k-2}$ are already constructed and
$Q_{i+1}$ is obtained by removing from the effective cycle
$(Q_i\circ D_{2,i+3-k})$ all irreducible components, not
containing the point $o$.\vspace{0.1cm}

Set $C=Q_{M-4}$. This is a 1-cycle, each component of which
contains the point $o$. We have $\mathop{\rm deg}C\leq
2^{k_2-2}\mathop{\rm deg}Q$, and for the multiplicity at the point
$o$ the equality
$$
\mathop{\rm mult}\nolimits_oC=2^{k-2}3^{k_2-2}\mathop{\rm
mult}\nolimits_oQ
$$
holds. Taking into account the choice of the subvariety $Q$, we
get
$$
\frac{\mathop{\rm mult}\nolimits_oC}{\mathop{\rm deg}C}>
\frac{2^{k-2}3^{k_2-2}}{2^{k_2-2}}
\cdot\frac{8}{d}=\frac{8}{9}.
$$
This estimate is obviously very strong, however it still does not
allow to get a contradiction: any line, passing through the point
$o$, satisfies this inequality. Now let us use the condition (R3),
bounding the input of the lines into the effective cycle
$C$.\vspace{0.3cm}

%%%%%%%%%%%%%%%%%%%%%%%%%%%%%%%%%%%%%%%%%%%%%%%%%%%%%%%%%%%%%%%%%%%%%
%%%%%%%%%%%%%%%%%%%%%%%%%%%%%%%%%%%%%%%%%%%%%%%%%%%%%%%%%%%%%%%%%%%%%
%%%%%%%%%%%%%%%%%%%%%%%%%%%%%%%%%%% subsection 1.3

{\bf 1.3. The effective 1-cycle, free from lines.} Now let us
write down $C=C_1+C_{\geq 2}$, where the support of the effective
1-cycle $C_1$ consists of lines, whereas the support of the
effective 1-cycle $C_{\geq 2}$ consists of the curves of degree
$\geq 2$. Obviously, $\mathop{\rm mult}\nolimits_oC_1=\mathop{\rm
deg}C_1$.\vspace{0.1cm}

{\bf Lemma 1.1.} {\it For any irreducible curve $\Gamma\subset V$
of degree $\mathop{\rm deg}\Gamma\geq 2$ the following inequality
holds:}
$$
\mathop{\rm mult}\nolimits_o\Gamma\leq\frac23\mathop{\rm deg}\Gamma.
$$
\vspace{0.1cm}

{\bf Proof.} By the condition (R1), the system of homogeneous
equations in the space ${\mathbb C}^{M+k}$
$$
q_{*,*}(z)=0,
$$
consisting of all irreducible components of all $k$ equations
$f_i$, defines either the origin $o\in{\mathbb C}^{M+k}$, or a
finite set of lines, passing through the point $o$. Obviously, a
line $t(a_1,\dots,a_{M+k})$ lies on $V$ if and only if
$q_{i,j}(a_*)=0$ for all $i,j$. Assume that the curve $\Gamma\ni
o$ is not a line and satisfies the inequality
$$
\mathop{\rm mult}\nolimits_o\Gamma>\frac23\mathop{\rm deg}\Gamma.
$$
It is clear that $\Gamma$ is contained in the support of any
tangent divisor $D\in\Sigma_1$ (that is, of such a divisor that
$D^+\in|H-2E|$) and hypertangent divisor $D\in\Sigma_2$ (that is,
$D^+\in|2H-3E|$). Therefore, on the curve $\Gamma$ the following
polynomials vanish identically:
$$
q_{i,1},\quad i=1,\dots,k,
$$
and
$$
q_{i,1}+q_{i,2},\quad i=k_1+1,\dots,k.
$$
Besides, since $\Gamma\subset V$, on that curve identically vanish
the polynomials
$$
f_i=q_{i,1}+q_{i,2},\quad i=1,\dots,k_1,
$$
and
$$
f_i=q_{i,1}+q_{i,2}+q_{i,3},\quad i=k_1+1,\dots,k.
$$
We conclude that all homogeneous polynomials $q_{i,j}$ vanish on
$\Gamma$. But then $\Gamma$ is a line. This contradiction proves
the lemma.\vspace{0.1cm}

{\bf Corollary 1.1.} {\it The following estimate holds:}
$$
\mathop{\rm mult}\nolimits_oC_{\geq 2}
\leq\frac23\mathop{\rm deg}C_{\geq 2}.
$$
\vspace{0.1cm}

To complete the proof of Theorem 3, it remains to estimate from
above the input of lines into the 1-cycle $C$, that is, the ratio
$\mathop{\rm deg}C_1/\mathop{\rm deg}C$.\vspace{0.1cm}

Let us consider generic tangent divisors $D_{1,i}$ that were used
to construct the curve $C$. If
$$
(\lambda_{1,i}q_{1,1}+\lambda_{2,i}q_{2,1}+\dots+\lambda_{k,i}q_{k,1})|_V
$$
is the equation of the divisor $D_{1,i}$, then obviously
$$
g_i=-(\lambda_{1,i}\bar{q}_{1,2}+\lambda_{2,i}\bar{q}_{2,2}+\dots+
\lambda_{k,i}\bar{q}_{k,2})\in\Sigma^{\mathbb T}_1
$$
is the equation of its projectivized tangent cone. Therefore,
$(g_1,\dots,g_{k-2})$ form a generic tuple, that is, generic
element of the space $(\Sigma^{\mathbb T}_1)^{\times (k-2)}$. In a
similar way, the projectivized tangent cone of the divisor
$D_{2,i}$, $i=1,\dots,k_2-2$, has the equation
$$
g_{k-2+i}=\sum^k_{j=1}\bar{s}_j\bar{q}_{j,2}+
\sum^{k_1+k_2}_{j=k_1+1}\lambda_j\bar{q}_{j,3}
\in\Sigma^{\mathbb T}_2
$$
(to simplify our notations, we omit in the right hand side of this
equality the second index for $\bar{s}_j$ and $\lambda_j$ that
indicates the dependence on the number $i$ of the divisor
$D_{2,i}$), whereas $(g_1,\dots,g_{M-4})\in{\cal P}$ is a generic
element of the space ${\cal P}$. Since at each step of the
construction of the curve $C$ the intersection with the
(hyper)tangent divisor $D_{i,j}$ is proper, the zero-dimensional
cycle $C_E=(C^+\circ E)$ is the scheme-theoretic intersection
$$
(Q_E\circ D_{1,1}\circ\dots\circ D_{1,k-2}
\circ D_{2,1}\circ\dots\circ D_{2,k_2-2}),
$$
that is, we get the equality of 0-cycles
$$
(C^+_1\circ E)+(C^+_{\geq 2}\circ E)=
\sum_{i\in I}m_i(R_i\circ D_{1,1}\circ\dots\circ D_{2,k_2-2}),
$$
in the right hand side we get the scheme-theoretic intersection
with all $M-4$ (hyper)tangent divisors $D_{i,j}$ that took part in
the procedure of constructing the 1-cycle $C$. For any line
$L\subset V$, passing through the point $o$, in the notations of
Sec. 0.2 we have
$$
(L^+\cap E)\in \Delta
$$
(the set $\Delta$ consists precisely of all tangent directions of
all lines on $V$, passing through the point $o$). Therefore, the
support of the 0-cycle $(C^+_1\circ E)$ is a subset of the finite
set $\Delta$. By the condition (R3), that is, by the estimate
(\ref{9april12.7}),
$$
\sum_{p\in\Delta}\sum_{i\in I}m_i\mathop{\rm dim}{\cal O}_{p,R_i}/
(g_1,\dots,g_{M-4})\leq 2^{k-4}3^{k_2-1}\mathop{\rm mult}\nolimits_oQ,
$$
since, obviously, the degree of the effective cycle $Q_E$ of
codimension three on $E$ is precisely $\mathop{\rm mult_o}Q$. The
more so,
\begin{equation}\label{9april12.4}
\delta=\mathop{\rm deg}(C^+_1\circ E)=\mathop{\rm deg}C_1
\leq 2^{k-4}3^{k_2-1}\mathop{\rm mult}\nolimits_oQ.
\end{equation}
Therefore, for the 1-cycle $C_{\geq 2}$ the equality
$$
\mathop{\rm mult}\nolimits_oC_{\geq 2}=
2^{k-2}3^{k_2-3}\mathop{\rm mult}\nolimits_oQ-\delta
$$
and the inequality
$$
\mathop{\rm deg}C_{\geq 2}\leq 2^{k_2-2}\mathop{\rm deg}Q-\delta
$$
hold. By easy computations from the inequality of Corollary 1.1 we
get the estimate
$$
\delta\geq 2^{k-2}3^{k_2-1}\mathop{\rm mult}\nolimits_oQ-
2^{k_2-1}\mathop{\rm deg}Q.
$$
Combining the last inequality with (\ref{9april12.4}), we obtain
finally
$$
2^{k_2-1}\mathop{\rm deg}Q\geq 2^{k-4}3^{k_2}\mathop{\rm mult}\nolimits_oQ.
$$
Recalling that $k=k_1+k_2$ and $d=2^{k_1}3^{k_2}$, we obtain from
there the inequality
$$
\mathop{\rm mult}\nolimits_oQ\leq\frac{8}{d}\mathop{\rm deg}Q,
$$
which contradicts the inequality (\ref{9april12.1}). This
contradiction excludes the third (and last) option $\mathop{\rm
codim}B\geq 4$ and completes the proof of Theorem 3.\vspace{0.1cm}

{\bf Remark 1.1.} In all previous papers (see, for instance,
\cite{Pukh07b,Pukh98b,Pukh03,Pukh04c}) the technique of
hypertangent divisors was used in a somewhat different way: at
each step when intersecting with a hypertangent divisor, we
selected an irreducibe component with the maximal ratio
$\mathop{\rm mult}\nolimits_o/\mathop{\rm deg}$, the other
components of the scheme-theoretic intersec\-tion were ignored. In
the argument given above, in order to estimate the input of lines
at the last step, we needed to control the whole cycle of
intersection, not only one of its components. This type of
arguments works for the previous problems, too: instead of
selecting the component with the maximal ratio $\mathop{\rm
mult}\nolimits_o/\mathop{\rm deg}$, one may consider the whole
cycle of the scheme-theoretic intersection, removing only the
components that do not contain the point $o$ (since the regularity
conditions ensure that the procedure of taking the intersection
with a hypertangent divisor is well defined in a neighborhood of
that point only).\vspace{0.5cm}

%%%%%%%%%%%%%%%%%%%%%%%%%%%%%%%%%%%%%%%%%%%%%%%%%%%%%%%%%%%%%%%%%%%%%%%%%%%%%%%%%%%%
%%%%%%%%%%%%%%%%%%%%%%%%%%%%%%%%%%%%%%%%%%%%%%%%%%%%%%%%%%%%%%%%%%%%%%%%%%%%%%%%%%%%
%%%%%%%%%%%%%%%%%%%%%%%%%%%%%%%%%%%%%%%%%%%%%%%%%%%%%%%%%%%%%%%%%%%%%%%%%%%%%%%%%%%%
%%%%%%%%%%%%%%%%%%%%%%%%%%%%%%%%%%%%%%%%%%%%%%%%%%%%%%%%%%%%%%%%%%%%%%%%%%%%%%%%%%%%
%%%%%%%%%%%%%%%%%%%%%%%%%%%%%%%%%%% SECTION 2

\begin{center}
\bf \S 2. The local multiplicities. I.\\ The spaces of tuples of
polynomials
\end{center}
\vspace{0.3cm}

In this section we give a precise formulation of the problem of
estimating the local multiplicity: introduce the space of tuples
of polynomials, define the local effective multiplicities of
intersection and formulate the main result (Theorem 4). The theory
developed in this and two subsequent sections is independent of \S
1 and is self-contained. The notations are also independent of \S
1.\vspace{0.3cm}

{\bf 2.1. The tuples of polynomials and effective multiplicities.}
Let us fix the complex coordinate space ${\mathbb
C}^{N}_{(z_1,\dots,z_N)}$, $N\geq 1$, which we consider as
embedded in the projective space ${\mathbb
P}^N_{(x_0:x_1:\dots:x_N)}$ as a standard affine chart $\{x_0\neq
0\}$, т.е. $z_i=x_i/x_0$. By the symbol ${\cal P}_{d,N}$ we denote
the space of homogeneous polynomials of degree $d\geq 1$ in the
variables $z_*$. For $e\leq d$ set
$$
{\cal P}_{[e,d],N}=\mathop{\oplus}\limits^d_{i=e}{\cal P}_{i,N},
$$
for instance, ${\cal P}_{[1,d],N}$ is the space of polynomials of
degree $\leq d$ without the constant term. On each of these spaces
we have an action of the matrix group $GL_N({\mathbb C})$ of
linear changes of coordinates. Let
$$
{\cal P}^{(n_1,n_2)}_N={\cal P}^{n_1}_{[1,2],N}\times{\cal
P}^{n_2}_{[1,3],N}=\prod^{n_1}_{i=1}{\cal P}^{}_{[1,2],N}
\times\prod^{n_2}_{i=1}{\cal P}_{[1,3],N}
$$
be the space of tuples
$(f_1,\dots,f_{n_1},f_{n_1+1},\dots,f_{n_1+n_2})$ of polynomials,
where the first $n_1$ polynomialsare of degree $\leq 2$ and the
subsequent $n_2$ polynomials are of degree at most 3. All
polynomials vanish at the point $o=(0,\dots,0)$. We assume that
the inequality $n_1+n_2=N$ holds.\vspace{0.1cm}

For an effective cycle $R$ of pure codimension $l$ on ${\mathbb
P}^N$ and a tuple of polynomials $(f_1,\dots,f_N)\in {\cal
P}^{(n_1,n_2)}_N$ define the {\it effective multiplicity}
$$
\mu_R(f_1,\dots,f_N)\in{\mathbb Z}_+\cup\{\infty\}
$$
in the following way. If $o\not\in\mathop{\rm Supp}R$, then we set
$\mu_R(f_*)=0$. If the closed set
$$
Z_R(f_1,\dots,f_N)=\{f_1=\dots=f_N=0\}\cap R
$$
is of positive dimension at the point $o$, then we set
$\mu_R(f_*)=\infty$. If none of these two cases takes place, then
$\mu_R(f_*)$ is a positive integer, to define which we need some
additional constructions. Let
$$
\Sigma_1(f_*)=<\!f_1,\dots,f_{n_1}\!>\subset{\cal P}_{[1,2],N}
$$
be the linear span of the quadratic polynomials. Let
$$
\Sigma_{12}(f_*)=
<f_1{\cal P}_{[0,1],N},\dots,f_{n_1}{\cal P}_{[0,1],N}>
\subset{\cal P}_{[1,3],N}
$$
be the linear space of all linear combinations of the polynomials
$f_1,\dots,f_{n_1}$ with polynomials of degree at most 1 in $z_*$
as coefficients. In particular, we have the following inclusion:
$\Sigma_1\subset\Sigma_{12}$. Let
$$
\Sigma_{22}(f_*)=<f_{n_1+1},\dots,f_{n_1+n_2}>\subset{\cal P}_{[1,3],N}
$$
be the linear span of the cubic polynomials and
$$
\Sigma_2(f_*)=\Sigma_{12}(f_*)+\Sigma_{22}(f_*)\subset{\cal P}_{[1,3],N}
$$
the sum of these two subspaces. Define the {\it polynomial span}
$$
[f_1,\dots,f_N]\subset{\cal P}^{(n_1,n_2)}_N
$$
of the tuple $(f_*)$ as the set of all such tuples
$(f^{\sharp}_1,\dots,f^{\sharp}_N)$ that
$f^{\sharp}_1,\dots,f^{\sharp}_{n_1}\in\Sigma_1(f_*)$ and
$f^{\sharp}_{n_1+1},\dots,f^{\sharp}_{n_1+n_2}\in\Sigma_2(f_*)$.
Obviously, $[f_*]$ is a closed irreducible subset of the space
${\cal P}^{(n_1,n_2)}_N$.\vspace{0.1cm}

Now for an irreducible subvariety $R\subset{\mathbb P}^N$ of
codimension $l$, $R\ni o$, set
$$
\mu_R(f_1,\dots,f_N)=
\mathop{\rm dim}{\cal O}_{o,R}/(f^{\sharp}_{l+1},\dots,f^{\sharp}_N)
$$
where $(f^{\sharp}_1,\dots,f^{\sharp}_N)\in[f_1,\dots,f_N]$ is a
generic tuple.For an arbitrary effective cycle $R=\Sigma_{j\in
J}r_jR_j$ of pure codimension $l$ we define $\mu_R(f_*)$ by
linearity, setting
$$
\mu_R(f_*)=\sum\limits_{j\in J}r_j\mu_{R_j}(f_*).
$$
The definition given above is equivalent to the following one
which works for any effective cycle of codimension $l$:
$$
\mu_R(f_*)=\mathop{\rm mult}\nolimits_o(\{f^{\sharp}_{l+1}=0\}\circ
\dots\circ\{f^{\sharp}_N=0\}\circ R),
$$
in the brackets it is the effective zero-dimensional cycle of the
scheme-theoretic intersection of the hypersurfaces
$\{f^{\sharp}_i=0\}$, where $i=l+1,\dots,N$, and the cycle $R$ in
a neighborhood of the point $o$.\vspace{0.1cm}

{\bf Remark 2.1.} To define the effective multiplicity, one needs
to form the polynomial span $[f_*]$, as certain components of the
cycle $R$ and the projectivized tangent cone to $R$ at the point
$o$ can be entirely contained in the divisors $\{f_i=0\}$. For
that reason the numbers of the form
$$
\mathop{\rm dim}{\cal O}_{o,R}/(f_{i_1},\dots,f_{i_{N-l}})
$$
can turn out to be strictly higher than the effective
multiplicity, correctly defined above, even if $R$ is an
irreducible subvariety, not contained in the hypersurfaces
$\{f_i=0\}$.\vspace{0.3cm}

%%%%%%%%%%%%%%%%%%%%%%%%%%%%%%%%%%%%%%%%%%%%%%%%%%%%%%%%%%%%%%%%%%%%%%%%%%%%%%%%%%%%
%%%%%%%%%%%%%%%%%%%%%%%%%%%%%%%%%%%%%%%%%%%%%%%%%%%%%%%%%%%%%%%%%%%%%%%%%%%%%%%%%%%%
%%%%%%%%%%%%%%%%%%%%%%%%%%%%%%%%%%% subsection 2.2

{\bf 2.2. The Chow varieties and the local effective
multiplicities.} By the symbol ${\mathbb H}_{l,N}(d)$ we denote
the Chow variety, parametrizing effective cycles of pure
codimension $l$ and degree $d$ on ${\mathbb P}^N$. Consider the
sets
$$
{\cal X}_{l,N}(m,d)\subset{\cal P}^{n_1}_{[1,2],N}
\times{\cal P}^{n_2}_{[1,3],N}\times{\mathbb H}_{l,N}(d),
$$
consisting of such tuples $((f_1,\dots,f_N),R)$, that
$$
\mu_R(f_1,\dots,f_N)\geq m\in{\mathbb Z}_+\cup\{\infty\}.
$$
It is easy to see that ${\cal X}_{l,N}(m,d)$ is a closed algebraic
set. Denoting by the symbol $\pi_{\cal P}$ the projection
$$
((f_1,\dots,f_N),R)\mapsto(f_1,\dots,f_N),
$$
and taking into account that the Chow varieties are projective, we
get that
$$
X_{l,N}(m,d)=\pi_{\cal P}({\cal X}_{l,N}(m,d))
\subset{\cal P}^{n_1}_{[1,2],N}\times{\cal P}^{n_2}_{[1,3],N}
$$
is a closed algebraic set. Explicitly, it consists of such tuples
$(f_1,\dots,f_N)$, for which there exists an effective cycle
$R\in{\mathbb H}_{l,N}(d)$, satisfying the inequality
$\mu_R(f_1,\dots,f_N)\geq m$.\vspace{0.1cm}

Let $B\subset{\cal P}^{n_1}_{[1,2],N}\times{\cal
P}^{n_2}_{[1,3],N}$ be an irreducible subvariety. We define the
local effective multiplicity, setting
$$
\mu^{\rm local}_{l,N}(B,d)=\mu_{l,N}(B,d)=\mathop{\rm
max}\limits_{m\in{\mathbb Z}_+ \cup\{\infty\}}\{m\,|\,B\subset
X_{l,N}(m,d)\}.
$$
The explicit meaning of this definition is as follows:
$\mu_{l,N}(B,d)=m$, if for a generic tuple $(f_*)\in B$ and any
effective cycle $R\in{\mathbb H}_{l,N}(d)$ the inequality
$$
\mu_R(f_1,\dots,f_N)\leq m
$$
holds and for at least one cycle $R\in{\mathbb H}_{l,N}(d)$ this
inequality turns into the equality. It is clear that if
$\mu_{l,N}(B,d)=\infty$, then for a generic (and thus for any)
tuple $(f_1,\dots,f_N)$ the set of its zeros $Z(f_1,\dots,f_N)$
has a component of positive dimension, passing through the point
$o$. The converse is also true: if there is such a component, for
$R$ one can take a subvariety such that $\mathop{\rm dim} (R\cap
Z(f_*))\geq 1$. \vspace{0.1cm}

{\bf Proposition 2.1.} {\it For any $d\geq 1$ the equality
$$
\mathop{\rm codim}X_{l,N}(\infty,d)=n_1+2n_2+1.
$$
holds.} \vspace{0.1cm}

{\bf Proof} is given in \S 5.\vspace{0.1cm}

Since in \S\S 2-4 only local multiplicities are considered, to
simplify the notations we omit the indication that they are local:
we write $\mu$ instead of $\mu^{\rm local}$ or $\mu_{\rm local}$;
in \S 5 we will define the global multiplicities $\mu_{\rm total}$
and for the local multiplicities we will use the notation
$\mu^{\rm local}$ or $\mu_{\rm local}$.\vspace{0.1cm}

Finally, set $\mu_{l,N}(a,d)=m$ for $a\in{\mathbb Z}_+$, if there
is an irreducible subvariety $B\subset{\cal P}^{(n_1,n_2)}_N$ of
codimension at most $a$, such that the equality $\mu_{l,N}(B,d)=m$
holds and the inequality
$$
\mathop{\rm codim}X_{l,N}(m+1,d)\geq a+1
$$
is satisfied.\vspace{0.1cm}

{\bf Remark 2.2.} By the definitions given above, if the
inequality $\mu_{l,N}(a,d)\leq m$ is satisfied, it means that
$\mathop{\rm codim} X_{l,N}(m+1,d)\geq a+1$. Explicitly: the set
of such tuples $(f_*)\in {\cal P}^{(n_1,n_2)}_N$, that there
exists an effective cycle $R$ of codimension $l$ and degree $d$,
satisfying the inequality $\mu_R(f_1,\dots, f_N)\geq m+1$, is of
codimension at least $a+1$ in the space ${\cal
P}^{(n_1,n_2)}_N$.\vspace{0.1cm}

Obviously, by Proposition 2.1 for any $a\geq n_1+2n_2+1$ we have
$\mu_{l,N}(a,d)=\infty$ and for $a\leq n_1+2n_2$ we have
$\mu_{l,N}(a,d)<\infty$. Starting from this moment we assume that
$a\leq n_1+2n_2$. (In the sequel, we will be interested mainly in
the case $a=n_1+n_2=N$.)\vspace{0.1cm}

By construction, the sets $X_{l,N}(m,d)$ are invariant under the
linear changes of coordinates (the group $GL_N({\mathbb C})$).
Besides, these sets are invariant under the action of another
group $G(n_1,n_2)$ on the space ${\cal P}^{(n_1,n_2)}_N$, which we
will now define. The group $G(n_1,n_2)$ is an extension
$$
\begin{array}{ccc}1\to & {\cal P}^{\times (n_1n_2)}_{[0,1],N} &
\to G(n_1,n_2)\to GL_{n_1}({\mathbb C})\times GL_{n_2}({\mathbb C})\to 1.\\
& & \\
& \parallel & \\  &  &  \\
& {\mathbb C}^{n_1n_2(N+1)} &
\end{array}
$$
More precisely, to an element $g\in G(n_1,n_2)$ corresponds a
triple of matrices $(A_{11},A_{12},A_{22})$, where
$$
A_{11}\in GL_{n_1}({\mathbb C}),\quad A_{22}\in GL_{n_2}({\mathbb
C}),\quad A_{12}\in\mathop{\rm Mat}\nolimits_{(n_1,n_2)}({\cal
P}_{[0,1],N}).
$$
Setting $A_{11}=\|\alpha_{ij}\|_{1\leq i,j\leq n_1}$,
$A_{22}=\|\beta_{ij}\|_{1\leq i,j\leq n_2}$ and
$$
A_{12}=\|\gamma_{ij}(z_1,\dots,z_N)\|_{1\leq i\leq n_1,1\leq j\leq n_2},
$$
we define the action
$$
g\colon f=(f_1,\dots,f_N)\mapsto f^g=(f^g_1,\dots,f^g_N)\in{\cal P}^{(n_1,n_2)}_N
$$
by the formulas
$$
f^g_i=\sum^{n_1}_{j=1}\alpha_{ij}f_j,\quad i=1,\dots,n_1,
$$
$$
f^g_{n_1+i}=
\sum^{n_2}_{j=1}\beta_{ij}f_{n_1+j}+\sum^{n_1}_{i=1}\gamma_{ij}(z_*)f_i,
\quad j=1,\dots,n_2.
$$
The closed subset $B\subset{\cal P}^{(n_1,n_2)}_N$ is said to be
{\it bi-invariant}, if it is invariant under the action of both
groups, $GL_N({\mathbb C})$ (changes of coordinates) and
$G(n_1,n_2)$. In particular, the sets $X_{l,N}(m,d$) are
bi-invariant (the multiplicities $\mu_R(f_*)$ are obviously
invariant under the action of both groups). Note that the group
$G(n_1,n_2)$ contains the subgroup $\overline{G}(n_1,n_2)\subset
GL_N({\mathbb C})$, corresponding to the tuples of matrices
$(A_{11},A_{12},A_{22})$ с
$$
A_{12}\in\mathop{\rm Mat}\nolimits_{(n_1,n_2)}({\mathbb C})\subset
\mathop{\rm Mat}\nolimits_{(n_1,n_2)}({\cal P}_{[0,1],N}).
$$
\vspace{0.3cm}

%%%%%%%%%%%%%%%%%%%%%%%%%%%%%%%%%%%%%%%%%%%%%%%%%%%%%%%%%%%%%%%%%%%%%%%%%%%%%%%%%%%%
%%%%%%%%%%%%%%%%%%%%%%%%%%%%%%%%%%%%%%%%%%%%%%%%%%%%%%%%%%%%%%%%%%%%%%%%%%%%%%%%%%%%
%%%%%%%%%%%%%%%%%%%%%%%%%%%%%%%%%%% subsection 2.3

{\bf 2.3. Reducing to the standard form.} For an irreducible
subvariety $B\subset{\cal P}^{(n_1,n_2)}_N$ we define its {\it
type}
$$
\tau(B)=((a_1,b_1),(a_2,b_2))\in{\mathbb Z}^{\times 4}_+,
$$
setting: $a_1=\mathop{\rm codim}(\mathop{\rm pr}_1(B)\subset{\cal
P}^{n_1}_{[1,2],N})$, where $\mathop{\rm pr}_1\colon{\cal
P}^{(n_1,n_2)}_N\to{\cal P}^{n_1}_{[1,2],N}$ is the projection
onto the first $n_1$ direct factors; $a_2=a-a_1$, where
$a=\mathop{\rm codim}(B\subset{\cal P}^{(n_1,n_2)}_N)$;
$$
b_1=n_1-\mathop{\rm rk}(df_1(o),\dots,df_{n_1}(o))
$$
for a tuple of general position $(f_1,\dots,f_{n_1})\in
\mathop{\rm pr}_1(B)$; $b_2=b-b_1$, where
$$
b=\varepsilon(B)=N-\mathop{\rm rk}(df_1(o),\dots,df_N(o))
$$
for a tuple of general position $(f_1,\dots,f_N)\in B$. If the
subvariety $B$ is $G(n_1,n_2)$(or at least
$\overline{G}(n_1,n_2)$)-invariant, then in a generic tuple of $N$
polynomials $(f_1,\dots,f_N)\in B$ the linear forms
\begin{equation}\label{oct11.1}
df_1(o),\dots,df_{n_1-b_1}(o),df_{n_1+1}(o)\dots,df_{N-b_2}(o)
\end{equation}
are linearly independent, the forms
$$
df_{n_1-b_1+1}(o),\dots,df_{n_1}(o)
$$
are linear combinations of $df_1(o),\dots,df_{n_1-b_1}(o)$, and
the forms
$$
df_{N-b_2+1}(o),\dots,df_N(o)
$$
are linear combinations of the forms (\ref{oct11.1}). Therefore,
for an arbitrary $G(n_1,n_2)$(or
$\overline{G}(n_1,n_2)$)-invariant irreducible subvariety $B$
there exists a non-empty Zarisky open subset $B^o\subset B$ {\it},
on which the following map of {\it reducing to the standard form}
is well defined:
\begin{equation}\label{oct11.2}
\rho\colon B^o\to{\cal P}^{n_1-b_1}_{[1,2],N}\times
{\cal P}^{b_1}_{2,N}\times{\cal P}^{n_2-b_2}_{[1,3],N}
\times{\cal P}^{b_2}_{[2,3],N},
\end{equation}
so that this map transforms a generic tuple $(f_1,\dots,f_N)$ into
the tuple of polynomials
$$
(f_1,\dots,f_{n_1-b_1},f^+_{n_1-b_1+1},
\dots,f^+_{n_1},f_{n_1+1},\dots,f_{N-b_2},f^+_{N-b_2+1},
\dots,f^+_N),
$$
where $f^+_i$ for $i\in\{n_1-b_1+1,\dots,n_1\}$ are obtained by
subtracting from $f_i$ the uniquely determined linear combination
of polynomials $f_1,\dots,f_{n_1-b_1}$, as a result of which
$df^+_i(o)=0$, and $f^+_i$ for $i\in\{N-b_2+1,\dots,N\}$ are
obtained by subtracting from $f_i$ the uniquely determined linear
combination of polynomials $f_j$,
$j\in\{1,\dots,n_1-b_1,n_1+1,\dots,N-b_2\}$, so that
$df^+_i(o)=0$. We denote the closure of the image $\rho(B^o)$ by
the symbol $\overline{B}$. If the subvariety $B$ is invariant with
respect to the action of $\overline{G}(n_1,n_2)$, then the
coefficients of the above mentioned linear combinations can be
arbitrary and for that reason
$$
\begin{array}{cl}
\mathop{\rm dim}B & =\mathop{\rm dim}\overline{B}+b_1(n_1-b_1)+b_2(N-b)= \\
  & =\mathop{\rm dim}\overline{B}+b(n_1-b_1)+b_2(n_2-b_2).
\end{array}
$$
On the other hand, the direct product in (\ref{oct11.2}) is of
codimension $bN$ in the space ${\cal P}^{(n_1,n_2)}_N$. Therefore,
the equality
\begin{equation}\label{oct11.3}
\mathop{\rm codim}\overline{B}=\mathop{\rm
codim}B-b_1(n_2+b_1)-b_2b
\end{equation}
holds. In the last formula the codimension of each of the
subvarieties is meant with respect to the {\it relevant} ambient
space: for $B$ it is ${\cal P}^{(n_1,n_2)}_N$, for $\overline{B}$
it is the direct product in (\ref{oct11.2}). Starting from this
moment, unless otherwise specified, the codimension is always
meant with respect to the natural ambient space. Sometimes, for
the convenience of the reader, we remind, with respect to what
ambient space the codimension is meant.\vspace{0.1cm}

Obviously, every fibre of the map $\rho\colon B^o\to \rho(B^o)$ is
${\mathbb C}^{(b_1(n_1-b_1)+b_2(N-b))}$.\vspace{0.1cm}

By construction,
$$
\mu_{l,N}(a,d)=\mathop{\rm max}\limits_B\mu_{l,N}(B,d),
$$
where the maximum is taken over all bi-invariant irreducible
subvarieties $B\subset{\cal P}^{(n_1,n_2)}_N$ of codimension at
most $a$. Our method of estimating the numbers $\mu_{l,N}(B,d)$
(and thus the numbers $\mu_{l,N}(a,d)$) is based on controlling
the type $\tau(B)$ of these subvarieties. It is easy to see that
the conditions
$$
\mathop{\rm rk}(df_1(o),\dots,df_{n_1}(o))\leq n_1-b_1
$$
and
$$
\mathop{\rm rk}(df_1(o),\dots,df_N(o))\leq N-b
$$
define in ${\cal P}^{(n_1,n_2)}_N$ a bi-invariant irreducible
subvariety of the type $((a^*_1,b_1),(a^*_2,b_2))$, where
$$
a^*_1=b_1(N+b_1-n_1)\quad\mbox{and}\quad a^*_2=b^2-a^*_1.
$$
Therefore, a bi-invariant irreducible subvariety $B$ of the type
$\tau(B)=((a_1,b_1),(a_2,b_2))$ can exist only if the inequalities
$$
a_1\geq b_1(N+b_1-n_1)\quad\mbox{and}\quad a=a_1+a_2\geq b^2
$$
hold.\vspace{0.1cm}

{\bf Proposition 2.2.} {\it Assume that a (bi-invariant)
subvariety $B$ is of the type $\tau(B)=((a_1,0),(a_2,0))$, that
is, $b=0$. Then
$$
\mu_{l,N}(B,d)=d
$$
for any $d\geq 1$.}\vspace{0.1cm}

{\bf Proof.} Let $R\ni o$ be an irreducible subvariety in
${\mathbb P}^N$, $(f_1,\dots,f_N)\in B$ a generic tuple of
polynomials. By assumption,
$$
\mathop{\rm rk}(df_1(o),\dots,df_N(o))=N.
$$
Set $\Pi=\{df_1(o)=\dots=df_{n_1}(o)=0\}\subset{\mathbb C}^N$,
this is a linear subspace of codimension $n_1$. Let
$$
(f^{\sharp}_1,\dots,f^{\sharp}_N)\in[f_1,\dots,f_N]
\subset{\cal P}^{(n_1,n_2)}_N
$$
be a generic element of the polynomial span of the tuple $(f_*)$.
By construction, the forms
$df^{\sharp}_1(o),\dots,df^{\sharp}_{n_1}(o)$ are linear forms of
general position, vanishing on the space $\Pi$, and
$df^{\sharp}_{n_1+1}(o),\dots,df^{\sharp}_N(o)$ are arbitrary
linear forms of general position.\vspace{0.1cm}

If $n_1\leq l$, then the inequality
$$
\mu_R(f_*)=\mathop{\rm dim}{\cal
O}_{o,R}/(f^{\sharp}_{l+1},\dots,f^{\sharp}_N) =\mathop{\rm
mult}\nolimits_oR
$$
holds, since the differentials $df^{\sharp}_i(o)$,
$i=l+1,\dots,N$, are forms of general position. If $n_1>l$, then
for the first $n_1-l$ differentials $df^{\sharp}_i(o)$, where
$i=l+1,\dots,n_1$, there is a (unique) constraint
$$
df^{\sharp}_i(o)|_{\Pi}\equiv 0.
$$
Let $T_oR\subset{\mathbb C}^N$ be the algebraic tangent cone to
$R$ at the point $o$. Each of its components has codimension $l$.
Since $\mathop{\rm codim}(\Pi\subset{\mathbb C}^N)=n_1$, for a
generic $f^{\sharp}_{l+1}$ one may assume that the linear form
$df^{\sharp}_{l+1}$ vanishes on no components of the cone $T_oR$.
Therefore,
$$
\mathop{\rm mult}\nolimits_o(R\circ\{f^{\sharp}_{l+1}=0\})=
\mathop{\rm mult}\nolimits_oR.
$$
For the same reasons, for generic
$f^{\sharp}_{l+1},\dots,f^{\sharp}_{n_1}$ we have the equality
$$
\mathop{\rm mult}\nolimits_o(R\circ\{f^{\sharp}_{l+1}=0\}\circ
\dots\circ\{f^{\sharp}_{l+1}=0\})=\mathop{\rm mult}\nolimits_oR.
$$
Now the equality $\mu_R(f_*)=\mathop{\rm mult}\nolimits_oR$ is
obtained as in the case $n_1\leq l$. Since $\mathop{\rm
mult}\nolimits_oR\leq\mathop{\rm deg}R$ and the equality is
attained, Proposition 2 is proved. Q.E.D.\vspace{0.1cm}

Note that if the multiplicity $\mu_{l,N}(B,d)$ is finite, then it
does not exceed the number $2^{a_1-l}3^{a_2}d$ for $a_1\geq l$ and
the number $3^{a-l}d$ for $a_1< l$. However, these estimates are
too weak for our purposes. Below the following key facts will be
shown.\vspace{0.1cm}

{\bf Theorem 4.} {\it The following inequalities are true:
$$
\mu_{l,N}(B,d)\leq \frac{1}{2\pi b}
\left(\frac{2a-b(b-1)}{b^2}e^2\right)^bd
<\frac{1}{2\pi b}\left(\frac{2a}{b^2}e^2\right)^bd.
$$
and}
$$
\mu_{l,N}(a,d)\leq\frac{e^2}{2\pi[\sqrt{a}]}(2e^2)^{[\sqrt{a}]}d.
$$
\vspace{0.1cm}

{\bf Theorem 5.} {\it For $a\geq 17$ the following estimate
holds:}
$$
\mu_{l,N}(a,d)\leq\frac{e^2}{2\pi[\sqrt{a}]}\left(\frac53 e^2\right)^{[\sqrt{a}]}d.
$$
\vspace{0.1cm}

Therefore, for a fixed $l$ and $a\to\infty$ the effective
multiplicity has exponential growth of the order $\sqrt{a}$, and
not of the order $a$, as in the aprioric estimates
above.\vspace{0.5cm}

%%%%%%%%%%%%%%%%%%%%%%%%%%%%%%%%%%%%%%%%%%%%%%%%%%%%%%%%%%%%%%%%%%%%%%%%%%%%%%%%%%%%
%%%%%%%%%%%%%%%%%%%%%%%%%%%%%%%%%%%%%%%%%%%%%%%%%%%%%%%%%%%%%%%%%%%%%%%%%%%%%%%%%%%%
%%%%%%%%%%%%%%%%%%%%%%%%%%%%%%%%%%%%%%%%%%%%%%%%%%%%%%%%%%%%%%%%%%%%%%%%%%%%%%%%%%%%
%%%%%%%%%%%%%%%%%%%%%%%%%%%%%%%%%%%%%%%%%%%%%%%%%%%%%%%%%%%%%%%%%%%%%%%%%%%%%%%%%%%%
%%%%%%%%%%%%%%%%%%%%%%%%%%%%%%%%%%% SECTION 3

\begin{center}
\bf \S 3. Local multiplicities. II.\\ The inductive method of
estimating
\end{center}
\vspace{0.3cm}

In this section, we construct the key procedure of estimating the
local multiplicity in terms of the local multiplicities of the
truncated tuples of polynomials.\vspace{0.3cm}

{\bf 3.1. Splitting off a direct factor.} If $b_2\geq 1$, then
denote by the symbol $\pi_2$ the projection
$$
\begin{array}{ccc}
{\cal P}^{n_1-b_1}_{[1,2],N}\times{\cal P}^{b_1}_{2,N} &
\times & {\cal P}^{n_2-b_2}_{[1,3],N}\times{\cal P}^{b_2}_{[2,3],N}\\
& \downarrow &\\
{\cal P}^{n_1-b_1}_{[1,2],N}\times{\cal P}^{b_1}_{2,N} &
\times & {\cal P}^{n_2-b_2}_{[1,3],N}\times{\cal P}^{b_2-1}_{[2,3],N}
\end{array}
$$
along the last direct factor ${\cal P}_{[2,3],N}$. If $b_1\geq 1$,
then by the symbol $\pi_1$ denote the similar projection along the
last quadratic direct factor ${\cal P}_{2,N}$ onto ${\cal
P}^{n_1-b_1}_{[1,2],N}\times{\cal P}^{b_1-1}_{2,N}\times{\cal
P}^{n_2}_{[1,3],N}\times{\cal P}^{b_2}_{[2,3],N}$. For a closed
set $\overline{B}$, constructed above by reducing to the standard
form, denote by the symbol $[\overline{B}]^{(i)}_{N-1}$ the
closure of the set $\pi_i(\overline{B})$ in the corresponding
ambient space, $i=1,2$. Furthermore, let $\lambda_2$ and
$\lambda_1$ denote the projections, complementary to $\pi_2$ and
$\pi_1$, that is, the projections {\it onto} the last direct
factor and the last quadratic direct factor, respectively. That is
to say, $\pi_2\times \lambda_2$ and $\pi_1\times \lambda_1$ are
the identity map of the direct product ${\cal
P}^{n_1-b_1}_{[1,2],N}\times{\cal P}^{b_1}_{2,N}\times {\cal
P}^{n_2-b_2}_{[1,3],N}\times{\cal
P}^{b_2}_{[2,3],N}$.\vspace{0.1cm}

Note that since the set $B$ is bi-invariant, the closed set
$\overline{B}$ is invariant with respect to the group
$GL_N({\mathbb C})$ of the linear changes of coordinates and to
the subgroup $G(n_1,b_1;n_2,b_2)\subset G(n_1,n_2)$, the elements
of which correspond to such triples of matrices
$(A_{11},A_{12},A_{22})$, that $A_{11}\in GL_{n_1-b_1}({\mathbb
C})$ (and these matrices act on the polynomials
$f_1,\dots,f_{n_1-b_1}$, mapping $f_{n_1-b_1+1},\dots,f_{n_1}$ to
themselves), $A_{22}\in GL_{n_2-b_2}({\mathbb C})$ (these matrices
act on $f_{n_1+1},\dots,f_{n_1+n_2-b_2}$) and
$A_{12}\in\mathop{\rm Mat}_{(n_1,n_2)}({\cal P}_{[0,1],N})$, where
only the polynomials $\gamma_{ij}$ с $j=1,\dots,n_2-b_2$ have an
arbitrary constant term, for $j\geq n_2-b_2+1$ these polynomials
are homogeneous: $\gamma_{ij}\in{\cal P}_{1,N}$. This trivial
remark is used below for estimating the parameters of the variety
$B_h\subset{\cal P}^{(n_1,n_2-1)}_{N-1}$ without special
comments.\vspace{0.1cm}

To simplify notations, assume that $b_2\geq 1$, consider the
projection $\pi_2$ and write $[\overline{B}]_{N-1}$ instead of
$[\overline{B}]^{(2)}_{N-1}$. The modifications required in the
case $i=1$ are obvious and we will give only the final result. In
all detail we consider the case $i=2$.\vspace{0.1cm}

For a generic tuple $(f_1,\dots,f_{N-1})\in[\overline{B}]_{N-1}$
by the symbol
$$
[\overline{B}]^N=[\overline{B}]^N(f_1,\dots,f_{N-1})\subset{\cal P}_{[2,3],N}
$$
we denote the fibre of the projection
$\pi_2|_{\overline{B}}\colon\overline{B}\to[\overline{B}]_{N-1}$.
Obviously,
$$
\mathop{\rm codim}\overline{B}=
\mathop{\rm codim}[\overline{B}]_{N-1}+\mathop{\rm codim}[\overline{B}]^N
$$
(recall: the codimension is meant with respect to the {\it
relevant} ambient space, for instance, for $[\overline{B}]^N$ it
is the fibre of the projection $\pi_2$, that is, ${\cal
P}_{[2,3],N})$. Set $\gamma_N=\gamma_N(B)=\mathop{\rm
codim}[\overline{B}]^N$. Since $\mathop{\rm codim}B\leq a$, we
obtain (see (\ref{oct11.3})) the estimate
$$
\mathop{\rm codim}[\overline{B}]_{N-1}\leq a-b_1(n_2+b_1)-b_2b-\gamma_N.
$$
From this, in particular, follows the inequality
\begin{equation}\label{oct11.4}
0\leq\gamma_N\leq a-b_1(n_2+b_1)-b_2b.
\end{equation}
Let $h(z_*)\in{\cal P}_{1,N}$ be a non-zero linear form,
$H=\{h=0\}\subset{\mathbb C}^N$ the corresponding hyperplane. By
one and the same symbol $\pi_h$ we will denote the projection of
all spaces ${\cal P}_{[i,j],N}$ onto ${\cal P}_{[i,j],N-1}$ and
the corresponding projections for direct products of these spaces,
where it is meant that an isomorphism $H\cong{\mathbb C}^{N-1}$ is
fixed and
$$
\pi_h\colon f\mapsto f|_H
$$
is the restriction of the polynomial onto $H$.\vspace{0.1cm}

Let $B_h\subset{\cal P}^{(n_1,n_2-1)}_{N-1}$ be the smallest
bi-invariant (in the sense of the latter space, that is, invariant
under the action of the groups $GL_{N-1}({\mathbb C})$ and
$G(n_1,n_2-1)$) closed set, containing the set
$\pi_h([\overline{B}]_{N-1})$. Obviously, $B_h$ is a bi-invariant
irreducible subvariety. For a generic tuple
$(f_1,\dots,f_{N-1})\in[\overline{B}]_{N-1}$ consider the
following tuple of linear forms
$$
df_1(o)|_H,\dots,df_{n_1-b_1}(o)|_H,df_{n_1+1}(o)|_H,\dots,df_{N-b_2}(o)|_H.
$$
The rank of this tuple is either $N-b$, or $N-b-1$. More
precisely, there exists an integer-valued vector
$$
\alpha=(\alpha_1,\alpha_2)\in\{(0,1),(-1,1),(0,0)\}
$$
such that for the parameters $(b^h_1,b^h_2)$ of the variety $B_h$
the following equalities hold:
$$
b^h_1=b_1-\alpha_1,\quad b^h_2=b_2-\alpha_2.
$$
Let $\tau(B_h)=((a^h_1,b^h_1),(a^h_2,b^h_2))$ be the type of the
subvariety $B_h$. Furthermore, set $b^h=b^h_1+b^h_2$ and
$a^h=a^h_1+a^h_2$. Obviously, $a^h$ is the codimension of the
subvariety $B_h$ in the ambient space ${\cal
P}^{(n_1,n_2-1)}_{N-1}$.\vspace{0.1cm}

To realize the inductive procedure defined in Sec. 3.3 below, we
will need, apart from $B_h$, another family of subvarieties in the
space ${\cal P}^{(n_1,n_2-1)}_{N-1}$, which we will now introduce,
again for the case $i=2$, with the obvious modifications in the
case $i=1$. Take any polynomial $f\in
\lambda_2(\overline{B})\subset {\cal P}_{[2,3],N}$ and consider
the closed set
$$
\overline{B}(f)=\lambda^{-1}_2(f)\cap \overline{B}\subset {\cal
P}^{n_1-b_1}_{[1,2],N}\times{\cal P}^{b_1}_{2,N}\times {\cal
P}^{n_2-b_2}_{[1,3],N}\times{\cal P}^{b_2-1}_{[2,3],N}.
$$
We will consider the product in the right hand side of the last
inclusion as the natural ambient space for $\overline{B}(f)$; in
particular, the codimension of $\overline{B}(f)$ is meant with
respect to that space. By construction, $\overline{B}(f)\neq
\emptyset$ and
$$
\mathop{\rm codim} \overline{B}\geq \mathop{\rm codim}
\overline{B}(f)+\mathop{\rm codim}\lambda_2(\overline{B}),
$$
so that we obtain the estimate
\begin{equation}\label{14may12.1}
\mathop{\rm codim} \overline{B}(f)\leq
a-b_1(n_2+b_1)-b_2b-\gamma(f)
\end{equation}
for some $\gamma(f)\in{\mathbb Z}_+$, $\gamma(f)\geq \mathop{\rm
codim}\lambda_2(\overline{B})$. Once again, take a non-zero linear
form $h(z_*)\in{\cal P}_{1,N}$ with $H=\{h=0\}$ the corresponding
hyperplane, fix an isomorphism $H\cong{\mathbb C}^{N-1}$ and let
$\pi_h$ mean the same as above.\vspace{0.1cm}

Let $B_h(f)\subset{\cal P}^{(n_1,n_2-1)}_{N-1}$ be the smallest
bi-invariant (in the sense of the latter space, that is, invariant
under the action of the groups $GL_{N-1}({\mathbb C})$ and
$G(n_1,n_2-1)$) closed set, containing the set
$\pi_h(\overline{B}(f))$. Obviously, $B_h(f)\subset B_h$ (since
$\overline{B}(f)\subset [\overline{B}]_{N-1}$) is a bi-invariant
closed subset, which we without loss of generality can assume to
be an irreducible subvariety. We will put off a discussion of the
invariants of the subvarieties $B_h(f)$ until Sec. 3.3, because we
will use the polynomials $f$ of a special form there. Now let us
study the subvariety $B_h$ in more detail.\vspace{0.3cm}

%%%%%%%%%%%%%%%%%%%%%%%%%%%%%%%%%%%%%%%%%%%%%%%%%%%%%%%%%%%%%%%%%%%%%%%%%%%%%%%%%%%%
%%%%%%%%%%%%%%%%%%%%%%%%%%%%%%%%%%%%%%%%%%%%%%%%%%%%%%%%%%%%%%%%%%%%%%%%%%%%%%%%%%%%
%%%%%%%%%%%%%%%%%%%%%%%%%%%%%%%%%%% subsection 3.2

{\bf 3.2. Estimating the codimension of the subvariety $B^h$.} How
the type of the subvariety $B$ changes when we restrict onto a
hyperplane, is shown in the following\vspace{0.1cm}

{\bf Proposition 3.1.} (i) {\it In the case $\alpha=(0,1)$ the
estimates}
$$
a^h_1\leq a_1-b_1,\quad a^h\leq a-(2b-1)-\gamma_N
$$
{\it hold.}\vspace{0.1cm}

(ii) {\it In the case $\alpha=(-1,1)$ the estimates}
$$
a^h_1\leq a_1,\quad a^h\leq a-b-\gamma_N
$$
{\it hold.}\vspace{0.1cm}

(iii) {\it In the case $\alpha=(0,0)$ the estimates}
$$
a^h_1\leq a_1-b_1,\quad a^h\leq a-b_1-b-\gamma_N
$$
{\it hold.}\vspace{0.1cm}

{\bf Proof.} We will study in full detail the cases (i) and (ii),
in the case (iii) we only give the computations. So assume that
$\alpha=(0,1)$. This means that the forms $df_i(o)$,
$i\in\{1,\dots,n_1-b_1,n_1+1,\dots,N-b_2\}$, remain linearly
independent when restricted onto $H$. We get the estimate
$$
\mathop{\rm codim}\pi_h([\overline{B}]_{N-1})\leq
\mathop{\rm codim}[\overline{B}]_{N-1}\leq
a-b_1(n_2+b_1)-b_2b-\gamma_N.
$$
The codimension in the left hand side is the codimension with
respect to the ambient space
\begin{equation}\label{oct11.5}
{\cal P}^{n_1-b_1}_{[1,2],N-1}\times{\cal P}^{b_1}_{2,N-1}
\times{\cal P}^{n_2-b_2}_{[1,3],N-1}\times{\cal P}^{b_2-1}_{[2,3],N-1},
\end{equation}
and moreover, for a generic tuple
$(g_1,\dots,g_{N-1})\in\pi_h([\overline{B}]_{N-1})$ the linear
forms $dg_i(o)$, $i\in\{1,\dots,n_1-b_1,n_1+1,\dots,N-b_2\}$, are
linearly independent. Now let us apply the procedure, inverse to
the procedure of reducing to the standard form: the set $B_h$, due
to its bi-invariance, certainly contains all the tuples of the
form $(g^+_1,\dots,g^+_{N-1})$, where $g^+_i=g_i$ for $1\leq i\leq
n_1=b_1$,
$$
g^+_i=g_i+\sum^{n_1-b_1}_{j=1}\lambda_{ij}g_j
$$
for $i=n_1-b_1+1,\dots,n_1$, where the coefficients $\lambda_{ij}$
are arbitrary. Furthermore, $g^+_i=g_i$ for $n_1+1\leq i\leq
N-b_2$ and
$$
g^+_i=g_i+\sum^{n_1-b_1}_{j=1}\lambda_{ij}g_j+\sum^{N-b_2}_{j=n_1+1}\lambda_{ij}g_j,
$$
for $i=N-b_2+1,\dots,N-1$, the coefficients $\lambda_{ij}$ are
arbitrary. For that reason, the inequality
$$
\mathop{\rm dim}B_h\geq\mathop{\rm
dim}\pi_h([\overline{B}]_{N-1})+ b_1(n_1-b_1)+(b_2-1)(N-b)
$$
holds. On the other hand, the dimension of the space ${\cal
P}^{(n_1,n_2-1)}_{N-1}$ is higher than the dimension of the space
(\ref{oct11.5}) by $(b-1)(N-1)$. As a result we get:
$$
\begin{array}{cl}
a^h\leq  & a-b_1(n_2+b_1)-b_2b-\gamma_N- \\
  &  -b_1(n_1-b_1)-(b_2-1)(N-b)+(b-1)(N-1)=\\
  &  =a-(2b-1)-\gamma_N,
\end{array}
$$
as it was claimed. The estimate for $a^h_1$ is obtained by similar
(but simpler) arguments. The procedure of reducing to the standard
form is well defined on the quadratic components, so that for the
projection $\mathop{\rm pr_1}([\overline{B}]_{N-1})$ we get the
estimate
$$
\mathop{\rm codim}\mathop{\rm pr_1}([\overline{B}]_{N-1})\leq a_1-b_1(n_2+b_1).
$$
Note that for this closed set the ambient space is the direct
product
$$
{\cal P}^{n_1-b_1}_{[1,2],N}\times{\cal P}^{b_1}_{2,N}.
$$
Now, restricting onto the hyperplane $H$ and applying, as we did
it above, the procedure, inverse to the procedure of reducing to
the standard form, we obtain the required estimate:
$$
a^h\leq a_1-b_1(n_2+b_1)-b_1(n_1-b_1)+b_1(N-1)=a_1-b_1,
$$
as required.\vspace{0.1cm}

Let us consider the case (ii). Here for a generic tuple
$(f_1,\dots,f_{N-1})\in[\overline{B}]_{N-1}$ the rank of the
system of linear forms $(df_*(o))$ drops by one when we restrict
onto $H$, and already the forms
$$
df_1(o)|_H,\, \dots,\, df_{n_1-b_1}(o)|_H
$$
are linearly independent. By the $G(n_1,b_1;n_2,b_2)$-invariance
of the set $[\overline{B}]_{N-1}$ we may assume that for a generic
tuple $(g_1,\dots,g_{N-1})\in\pi_h([\overline{B}]_{N-1})$ the
first $n_1-b_1-1$ forms $dg_i(o)$, $i=1,\dots,n_1-b_1-1$, are
linearly independent, and $dg_{n_1-b_1}(o)$ is their linear
combination, and moreover, for any
$\lambda_1,\dots,\lambda_{n_1-b_1-1}$
$$
(g_1,\dots,g_{n_1-b_1-1},g_{n_1-b_1}+
\sum^{n_1-b_1-1}_{i=1}\lambda_ig_i,g_{n_1-b_1+1},\dots,g_{N-1})
\in\pi_h([\overline{B}]_{N-1}).
$$
Therefore, to the Zariski open subset
$$
\mathop{\rm rk}<dg_1(o),\dots,dg_{n_1-b_1-1}(o)>=n_1-b_1-1
$$
one can apply the map of reducing to the standard form in the
component $g_{n_1-b_1}$ and, taking the closure, obtain the
irreducible subset
$$
\overline{\pi_h([\overline{B}]_{N-1})}
\subset{\cal P}^{n_1-b_1-1}_{[1,2],N-1}\times
{\cal P}^{b_1+1}_{2,N-1}\times{\cal P}^{n_2-b_2}_{[1,3],N-1}\times
{\cal P}^{b_2-1}_{[2,3],N-1}
$$
of codimension
$$
\mathop{\rm codim}\pi_h([\overline{B}]_{N-1})-N+n_1-b_1.
$$
After that, one can apply the procedure, inverse to that of
reducing to the standard form, similar to how it was done in the
case (i), and obtain for the codimension of the set $B_h$ the
estimate
$$
\mathop{\rm codim}B_h\leq a-b_1(n_2+b_1)-b_2b-\gamma_N-N+n_1-b_1+b(N-1)-
$$
$$
-(b_1+1)(n_1-b_1-1)-(b_2-1)(N-b-1)=a-b-\gamma_N,
$$
as we claimed. For $a^h_1$ the arguments are similar, but simpler,
since we take into account only the quadratic polynomials. Again
one have to reduce the generic tuple $(g_1,\dots,g_{n_1})$ to the
standard form in the component $g_{n_1-b_1}$ and after that apply
to the resulting set the procedure, inverse to that of reducing to
the standard form. As a result, we obtain the inequality
$$
\begin{array}{cl}
a^h_1\leq  &  a_1-b_1(n_2+b_1)-N+n_1-b_1- \\
  &  -(b_1+1)(n_1-b_1-1)+(b_1+1)(N-1)=a_1,
\end{array}
$$
as it was claimed (recall that $N=n_1+n_2$).\vspace{0.1cm}

Finally, let us consider the case (iii). As in the previous case,
one has to complete reducing the set $\pi_h([\overline{B}]_{N-1})$
to the standard form, now in the component $g_{N-b_2}$, since for
a generic tuple $(g_1,\dots,g_{N-1})$ the linear forms
$$
dg_1(o),\dots,dg_{n_1-b_1}(o),dg_{n_1+1}(o),\dots,dg_{N-b_2-1}(o)
$$
are linearly independent, and $dg_{N-b_2}(o)$ is their linear
combination. After that we apply the procedure, inverse to that of
reducing to the standard form. We obtain the inequality:
$$
\begin{array}{cl}
a^h\leq  & a-b_1(n_2+b_1)-b_2b-\gamma_N- \\
  &  -(N-1)+(N-b-1)-  \\
  &  -b_1(n_1-b_1)-b_2(N-b-1)+b(N-1),
\end{array}
$$
where the second line corresponds to the reduction to the standard
form and the third to the inverse procedure. Simplifying, we
obtain the estimate
$$
a^h\leq a-b_1-b-\gamma_N,
$$
as it was claimed. For the codimension $a^h_1$ in the quadratic
components this case is very easy, as there is no need to reduce
to the standard form. Applying the procedure, inverse to the
procedure of reducing to the standard form, we obtain the
inequality:
$$
a^h_1\leq a_1-b_1(n_2+b_1)-b_1(n_1-b_1)+b_1(N-1)=a_1-b_1,
$$
as it was claimed. Q.E.D. for Proposition 3.1.\vspace{0.1cm}

Now let us consider the problem of estimating the codimension when
a quadratic factor is being split off. Here
$(b^h_1,b^h_2)=(b_1,b_2)-(\alpha_1,\alpha_2)$, where
$\alpha=(\alpha_1,\alpha_2)$ can take the values $(1,0),(1,-1)$
and $(0,0)$. The inequalities for $a^h$ and $a^h_1$ are obtained
by word for word the same arguments as when a cubic factor was
split off. We give the final result.\vspace{0.1cm}

{\bf Proposition 3.2.} (i) {\it In the case $\alpha=(1,0)$ the
estimates}
$$
a^h_1\leq a_1-n_2-2b_1+1-\gamma_N,\quad a^h\leq
a-n_2-2b_1-b_2+1-\gamma_N
$$
{\it hold.}\vspace{0.1cm}

(ii) {\it In the case $\alpha=(1,-1)$ the estimates}
$$
a^h_1\leq a_1-n_2-2b_1+1-\gamma_N,\quad a^h\leq
a-n_2-2b_1+1-\gamma_N
$$
{\it hold.}\vspace{0.1cm}

(iii) {\it In the case $\alpha=(0,0)$ the estimates}
$$
a^h_1\leq a_1-n_2-b_1-\gamma_N,\quad a^h\leq a-n_2-b_1-\gamma_N
$$
{\it hold.}\vspace{0.1cm}

{\bf Proof} is left to the reader.\vspace{0.1cm}

To conclude, we remind the reader that when a quadratic factor is
being split off, the ambient space for the subvariety $B^h$ is
${\cal P}^{(n_1-1,n_2)}_{N-1}$, and the codimension $a^h$ is the
codimension of the subvariety $B_h$ with respect to the latter
space.\vspace{0.3cm}

%%%%%%%%%%%%%%%%%%%%%%%%%%%%%%%%%%%%%%%%%%%%%%%%%%%%%%%%%%%%%%%%%%%
%%%%%%%%%%%%%%%%%%%%%%%%%%%%%%%%%%%%%%%%%%%%%%%%%%%%%%%%%%%%%%%%%%%
%%%%%%%%%%%%%%%%%%%%%%%%%%%%%%%%%%% subsection 3.3

{\bf 3.3. The main inductive estimate.} Let us come back to the
main problem of estimating the multiplicities $\mu_{l,N}(B,d)$ for
an irreducible bi-invariant subvariety $B\subset{\cal
P}^{(n_1,n_2)}_N$ of codimension $a\geq 1$. Let
$\tau(B)=((a_1,b_1),(a_2,b_2))$ be the type of the subvariety $B$,
where $b=b_1+b_2\geq 1$ (otherwise by Proposition 2.2 there is
nothing to estimate). The following key fact makes it possible to
estimate the multiplicity $\mu_{l,N}$ from above in terms of
similar multiplicities for ${\mathbb P}^{N-1}$, which gives an
inductive  (in $N$) procedure of estimating
multiplicities.\vspace{0.1cm}

{\bf Proposition 3.3.} (i) {\it Assume that $b_2\geq 1$. Then
there exist non-zero linear forms $h_1,h_2\in{\cal P}_{1,N}$,
depending on the variety $B$ only and a set of non-negative
integers
$$
\left(\begin{array}{cc}d_{11} & d_{12}\\
d_{21} & d_{22}\end{array}\right)\in
\mathop{\rm Mat}\nolimits_{2,2}({\mathbb Z}_+),
$$
satisfying the equalities $d_{11}+d_{12}=d_{21}+d_{22}=d$, such
that the subvarieties}
$$
B_1=B_{h_1}\subset{\cal P}^{(n_1,n_2-1)}_{N-1}\quad \mbox{\it
and}\quad B_2=B_{h_2}(h_1h_2)\subset{\cal P}^{(n_1,n_2-1)}_{N-1}
$$
{\it satisfy the inequality}
\begin{equation}\label{oct11.6}
\begin{array}{cl}
\mu_{l,N}(B,d)\leq  &  \mu_{l,N-1}(B_1,d_{11})+\mu_{l-1,N-1}(B_1,d_{12})+\\
                    & +\mu_{l,N-1}(B_2,d_{21})+\mu_{l-1,N-1}(B_2,d_{22}).
\end{array}
\end{equation}
{\it Moreover, the subvariety $B_1$ is of the type}
$$
\tau(B_1)=((a_{11},b_{11}),(a_{12},b_{12}))\quad \mbox{\it
with}\quad (b_{11},b_{12})=(b_1,b_2-1).
$$
\vspace{0.1cm}

(ii) {\it Assume that $b_1\geq 1$. Then there exist non-zero
linear forms $h_1,h_2$, depending on $B$ only, and a set of
non-negative integers $(d_{ij})_{1\leq i,j\leq 2}$, satisfying the
equalities $d_{11}+d_{12}=d_{21}+d_{22}=d$, such that for the
subvarieties
$$
B_1=B_{h_1}\subset{\cal P}^{(n_1-1,n_2)}_{N-1}\quad\mbox{\it
and}\quad B_2= B_{h_2}(h_1h_2)\subset{\cal P}^{(n_1-1,n_2)}_{N-1}
$$
{\it the inequality (\ref{oct11.6}) is satisfied. Moreover, the
subvariety $B_1$ has the type}
$\tau(B_1)=((a_{11},b_{11}),(a_{12},b_{12}))$ with
$(b_{11},b_{12})=(b_1-1,b_2)$.}\vspace{0.1cm}

{\bf Proof.} Let us consider the case (i) in full detail. The
second case is considered in an absolutely similar way. The
bi-invariant irreducible subvariety $B$ is fixed. Let $h_1\in{\cal
P}_{1,N}$ be a linear form of general position, where the
genericity is understood in the following sense: for a generic
tuple $(f_1,\dots,f_N)\in B$ the linear subspace
$$
\{df_1(o)=\dots=df_N(o)=0\}
$$
(which is of dimension $b=b_1+b_2$) is not contained in the
hyperplane $\{h_1=0\}$. Let
$$
\Pi=\{h_1(z_*)h(z_*)\,|\, h\in{\cal P}_{1,N}\}\subset{\cal P}_{2,N}
$$
be the linear space of reducible homogeneous quadratic
polynomials, divisible by $h_1$. Note that ${\cal
P}_{2,N}\subset{\cal P}_{[2,3],N}$, so that $\Pi$ can be
considered as a linear subspace in ${\cal P}_{[2,3],N}$.
Obviously, $\mathop{\rm dim}\Pi=N$. Set
$$
{\cal P}_{\Pi}={\cal P}^{n_1-b_1}_{[1,2],N}\times{\cal P}^{b_1}_{2,N}
\times{\cal P}^{n_2-b_2}_{[1,3],N}\times{\cal P}^{b_2-1}_{[2,3],N}\times\Pi.
$$
As we explained above, ${\cal P}_{\Pi}$ is a closed irreducible
subset of the space
$$
{\cal P}_{\Pi}={\cal P}^{n_1-b_1}_{[1,2],N}\times{\cal P}^{b_1}_{2,N}
\times{\cal P}^{n_2-b_2}_{[1,3],N}\times{\cal P}^{b_2}_{[2,3],N}.
$$
\vspace{0.1cm}

{\bf Lemma 3.1.} (i) {\it The intersection $\overline{B}\cap{\cal
P}_{\Pi}$ is non-empty and its codimension in ${\cal P}_{\Pi}$ is
not higher than $\mathop{\rm codim}\overline{B}$.}\vspace{0.1cm}

(ii) {\it The closure $\overline{\pi_2(\overline{B}\cap{\cal
P}_{\Pi})}$ coincides with $[\overline{B}]_{N-1}$.}\vspace{0.1cm}

{\bf Proof.} Let us show the claim (i). By the bi-invariance the
closed set $\overline{B}$ contains the zero tuple
$(0,\dots,0)\in{\cal P}^{(n_1,n_2)}_N$. Therefore
$\overline{B}\cap{\cal P}\neq\emptyset$, and the rest is
obvious.\vspace{0.1cm}

Now let us show the claim (ii). Recall that $\pi_2$ is the
projection along the last direct factor ${\cal P}_{[2,3],N}$. By
the equality (\ref{oct11.3}) the codimension of the set
$\overline{B}$ is strictly smaller than $N$. Therefore, for a
tuple of general position
$(f_1,\dots,f_{N-1})\in[\overline{B}]_{N-1}$ the intersection of
the fibre $[\overline{B}]^N(f_1,\dots,f_{N-1})$ with the subspace
$\Pi$ in the space  ${\cal P}_{[2,3],N}$ is non-empty and for that
reason has a positive dimension. (The non-emptiness is again a
consequence of being bi-invariant:
$$
(f_1,\dots,f_{N-1},0)\in\overline{B}.)
$$
This proves the claim (ii). Q.E.D. for the lemma.\vspace{0.1cm}

{\bf Remark 3.1.} Without loss of generality, we may assume that
the irreducible subvariety $B$ is an irreducible component of the
closed set $X_{l,N}(m,d)$.\vspace{0.1cm}

The equality $\mu_{l,N}(B,d)=m$ means that for a generic tuple
$(f_1,\dots,f_N)\in\overline{B}$ and any effective cycle $R$ of
codimension $l$ and degree $d$ the inequality
$\mu_R(f_1,\dots,f_N)\leq m$ holds, and moreover, for some cycle
$R$ (depending on the tuple $(f_*)$) the equality holds.
Therefore, for a generic tuple
$$
(f_1,\dots,f_{N-1},h_1(z_*)h_2(z_*))\in\overline{B}\cap{\cal P}_{\Pi}
$$
there is en effective cycle $R$, depending on that tuple, for
which the inequality
$$
\mu_R(f_1,\dots,f_{N-1},h_1h_2)\geq m
$$
is satisfied. The form $h_2(z_*)$, as we explained above, is
non-zero. Let $H_i=\{h_i=0\}$, $i=1,2,$ be the corresponding
hyperplanes. Define the effective cycles $R_{ij}$ of codimension
$l$, $i,j=1,2,$ by the conditions:
\begin{itemize}

\item $R=R_{11}+R_{12}=R_{21}+R_{22}$,

\item $\mathop{\rm Supp}R_{12}\subset H_1$ and none of the
irreducible components of the cycle $R_{11}$ is contained in
$H_1$,

\item $\mathop{\rm Supp}R_{22}\subset H_2$ and none of the
irreducible components of the cycle $R_{21}$ is contained in
$H_2$.
\end{itemize}

Set $d_{ij}=\mathop{\rm deg}R_{ij}$. Obviously,
$d=d_{11}+d_{12}=d_{21}+d_{22}$. By what was said above, the
inequality
\begin{equation}\label{oct11.7}
m\leq\mu_R(f_1,\dots,f_{N-1},h_1)+\mu_R(f_1,\dots,f_{N-1},h_2)
\end{equation}
holds, and for any $i\in\{1,2\}$
\begin{equation}\label{oct11.8}
\mu_R(f_1,\dots,f_{N-1},h_i)\leq\mu_{R_{i1}}(f_1,\dots,f_{N-1},h_i)+
\mu_{R_{i2}}(f_1,\dots,f_{N-1},h_i).
\end{equation}
Since the irreducible components of the cycle $R_{i1}$ are not
contained in the hyperplane $H_i$, and $\mathop{\rm
Supp}R_{i2}\subset H_i$, for the first and second multiplicities
in the right hand side of the last inequality we get the estimates
\begin{equation}\label{oct11.9}
\mu_{R_{i1}}(f_1,\dots,f_{N-1},h_i)\leq
\mu_{(R_{i1}\circ H_i)}(f_1|_{H_i},\dots,f_{N-1}|_{H_i})
\end{equation}
and
\begin{equation}\label{oct11.10}
\mu_{R_{i2}}(f_1,\dots,f_{N-1},h_i)\leq
\mu_{R_{i2}}(f_1|_{H_i},\dots,f_{N-1}|_{H_i}),
\end{equation}
respectively. Note that $(R_{i1}\circ H_i)$ is an effective cycle
of codimension $l$ and degree $d_{i1}$ on $H_i\cong{\mathbb
P}^{N-1}$, and $R_{i2}$ is an effective cycle of codimension $l-1$
on $H_i\cong{\mathbb P}^{N-1}$. Furthermore,
$$
(f_1|_{H_1},\dots,f_{N-1}|_{H_1})\in\pi_{h_1}([\overline{B}]_{N-1})
$$
and
$$
(f_1|_{H_2},\dots,f_{N-1}|_{H_2})\in\pi_{h_2}(\overline{B}(h_1h_2))
$$
are generic tuples, so that the inequalities (\ref{oct11.9}) and
(\ref{oct11.10}) remain true, if in the right hand side that tuple
is replaced by a generic tuple of polynomials
$$
(g_1,\dots,g_{N-1})\in B_i.
$$
Now the inequality (\ref{oct11.6}) is a direct corollary of
(\ref{oct11.7}), (\ref{oct11.8}), (\ref{oct11.9}) and
(\ref{oct11.10}).\vspace{0.1cm}

Finally, by the genericity of the form $h_1$ we have
$$
\mathop{\rm rk}(dg_1(o),\dots,dg_{N-1}(o))=N-b,
$$
as it was claimed. Q.E.D. for Proposition 3.3.\vspace{0.1cm}

In order to make our inductive procedure a working one, it remains
to estimate the codimension of the subvarieties $B_2$ and their
type. Again we consider more closely the case of the projection
onto the last factor (that is, $(f_*)\mapsto f_N$), corresponding
to the part (i) of Proposition 3.3. In the case of the projection
onto the last quadratic factor (that is, $(f_*)\mapsto f_{n_1}$)
the modifications are obvious and we give only the final result.
Set $\gamma=\gamma(h_1h_2)$ (see Sec. 3.1).\vspace{0.1cm}

{\bf Proposition 3.4.} {\it Let $h_2\in{\cal P}_{1,N}$ be a
generic linear form such that $h_1h_2\in\lambda_2(\overline{B}\cap
{\cal P}_{\Pi})$. There exists and integer-valued vector
$$
\alpha=(\alpha_1,\alpha_2)\in\{(0,1),(-1,1),(0,0)\}
$$
such that for the parameters $({\bar b}_1,{\bar b}_2)$ of the
subvariety $B_{h_2}(h_1h_2)$ the following equalities hold:
$$
{\bar b}_1=b_1-\alpha_1,\quad {\bar b}_2=b_2-\alpha_2.
$$
For the codimension ${\bar a}$ of the subvariety $B_{h_2}(h_1h_2)$
in the ambient space ${\cal P}^{(n_1,n_2-1)}_{N-1}$ the following
estimates hold.}\vspace{0.1cm}

(i) {\it In the case $\alpha=(0,1)$:} ${\bar a}\leq
a-(2b-1)-\gamma$.\vspace{0.1cm}

(ii) {\it In the case $\alpha=(-1,1)$} ${\bar a}\leq
a-b-\gamma$.\vspace{0.1cm}

(iii) {\it In the case $\alpha=(0,0)$:} ${\bar a}\leq
a-b_1-b-\gamma$.\vspace{0.1cm}

{\bf Proof} is almost word for word the same as that of
Proposition 3.1. Although the set $\overline{B}(h_1h_2)$ may be
not invariant under the linear changes of coordinates, it is still
invariant under the operations of taking linear combinations, in
the same way as the set $[\overline{B}]_{N-1}$. By Lemma 3.1,
$$
\overline{\mathop{\bigcup}\limits_{h_1h_2\in\lambda_2(\overline{B}\cap{\cal
P}_{\Pi})} \overline{B}(h_1h_2)}=[\overline{B}]_{N-1}
$$
and for this reason for a generic tuple
$(g_1,\dots,g_{N-1})\in\overline{B}(h_1h_2)$ the linear forms
$$
dg_1(o),\dots,dg_{n_1-b_1}(o), dg_{n_1+1}(o),\dots,dg_{N-b_2}(o)
$$
are linearly independent. Taking this into account, the proof of
Proposition 3.1 works word for word, given the inequality
(\ref{14may12.1}), and with simplifications as we claim nothing
about the parameter ${\bar a}_1$ of the full type of
$B_{h_2}(h_1h_2)$. Q.E.D.\vspace{0.1cm}

Now let us formulate the result for the case when a quadratic
factor is split off.\vspace{0.1cm}

{\bf Proposition 3.5.} {\it Let $h_2\in{\cal P}_{1,N}$ be a
generic linear form such that $h_1h_2\in\lambda_1(\overline{B}\cap
{\cal P}_{\Pi})$. There exists and integer-valued vector
$$
\alpha=(\alpha_1,\alpha_2)\in\{(1,0),(1,-1),(0,0)\}
$$
such that for the parameters $({\bar b}_1,{\bar b}_2)$ of the
subvariety $B_{h_2}(h_1h_2)$ the following equalities hold:
$$
{\bar b}_1=b_1-\alpha_1,\quad {\bar b}_2=b_2-\alpha_2.
$$
For the codimension ${\bar a}$ of the subvariety $B_{h_2}(h_1h_2)$
in the ambient space ${\cal P}^{(n_1-1,n_2)}_{N-1}$ the following
estimates hold.}\vspace{0.1cm}

(i) {\it In the case $\alpha=(1,0)$:} ${\bar a}\leq
a-n_2-2b_1-b_2+1-\gamma$.\vspace{0.1cm}

(ii) {\it In the case $\alpha=(1,-1)$} ${\bar a}\leq
a-n_2-2b_1+1-\gamma$.\vspace{0.1cm}

(iii) {\it In the case $\alpha=(0,0)$:} ${\bar a}\leq
a-n_2-b_1-\gamma$.\vspace{0.1cm}

{\bf Proof} is almost word for word the same as that of
Proposition 3.2, and follows the same procedure as was used in the
proof of Proposition 3.1. Q.E.D.\vspace{0.5cm}

%%%%%%%%%%%%%%%%%%%%%%%%%%%%%%%%%%%%%%%%%%%%%%%%%%%%%%%%%%%%%%%%%%%%%%%%%%%%%%%%%%%%
%%%%%%%%%%%%%%%%%%%%%%%%%%%%%%%%%%%%%%%%%%%%%%%%%%%%%%%%%%%%%%%%%%%%%%%%%%%%%%%%%%%%
%%%%%%%%%%%%%%%%%%%%%%%%%%%%%%%%%%%%%%%%%%%%%%%%%%%%%%%%%%%%%%%%%%%%%%%%%%%%%%%%%%%%
%%%%%%%%%%%%%%%%%%%%%%%%%%%%%%%%%%%%%%%%%%%%%%%%%%%%%%%%%%%%%%%%%%%%%%%%%%%%%%%%%%%%
%%%%%%%%%%%%%%%%%%%%%%%%%%%%%%%%%%% SECTION 4

\begin{center}
\bf \S 4. Local multiplicities. III.\\ Explicit estimates
\end{center}
\vspace{0.3cm}

In this section, using the inductive procedure, developed in \S 3,
we obtain explicit estimates for the local multiplicity. We
consider separately the cases of small values $b=1,2$ and small
codimensions $a\leq 36$. We prove Theorems 4 and 5.\vspace{0.3cm}

{\bf 4.1. An estimate for the multiplicity in the case
$b=b_1+b_2=1$.} We did see above that for an irreducible
subvariety $B\subset{\cal P}^{(n_1,n_2)}_N$ with $b=0$ the
equality $\mu_{l,N}(B,d)=d$ holds (Proposition 2.2). Let us
consider the case $b=1$, the next in complexity. Assume for
certainty that $b_2=1$, $b_1=0$. By Proposition 3.3, the
inequality
$$
\mu_{l,N}(B,d)\leq
d+\mu_{l,N-1}(B_2,d_{21})+\mu_{l-1,N-1}(B_2,d_{22})
$$
holds, since $(b_{11},b_{12})=(0,0)$. Let
$\tau(B_2)=((a_{21},b_{21}),(a_{22},b_{22}))$ be the type of the
variety $B_2$. If $b_{21}=b_{22}=0$, then we get the estimate
$\mu_{l,N}(B,d)\leq 2d$. Assume that $b_{21}+b_{22}=1$. Setting
$B_2=B^{(1)}$, let us apply Proposition 3.3 to that subvariety.
Iterating this construction, we obtain a chain of subvarieties
$$
B^{(1)},\dots,B^{(k)},
$$
of the type
$\tau(B^{(i)})=((a^{(i)}_1,b^{(i)}_1),(a^{(i)}_2,b^{(i)}_2))$ with
$b^{(i)}=b^{(i)}_1+b^{(i)}_2=1$. Here $B^{(i+1)}=B^{(i)}_2$ in the
sense of Proposition 3.3. The varieties $B^{(i)}$ are irreducible
bi-invariant subvarieties of the space ${\cal
P}^{(n^{(i)}_1,n^{(i)}_2)}_{N-i}$, the corresponding subvariety
$B^{(i)}_1$ has the type $((a^{(i)}_{11},0),(a^{(i)}_{12},0))$ and
its input into the estimate of the multiplicity $\mu_{l,N}$ is
known. After $k$ steps we get the inequality неравенство
$$
\mu_{l,N}(B,d)\leq  kd+
\sum^{{\rm min}\{k,l\}}_{j=0}\mu_{l-j,N-k}(B^{(k)},d_{k,j}),
$$
where $d_{k,0}+\dots+d_{k,{\rm min}\{k,l\}}=d$. We used the
obvious inequality
\begin{equation}\label{oct11.11}
\mu_{l,N}(B,d')+\mu_{l,N}(B,d'')\leq\mu_{l,N}(B,d'+d'').
\end{equation}

By Propositions 3.1 and 3.2 for the codimension
$a^{(i)}=\mathop{\rm codim}B^{(i)}$ we have the inequality
$a^{(i+1)}\leq a^{(i)}-1$, and for the variety $B^{(i)}$ to exist,
the inequality
$$
a^{(i)}-b^{(i)}_1(k^{(i)}_2+b^{(i)}_1)-b^{(i)}_2b^{(i)}\geq 0
$$
should be satisfied. Therefore, $a\geq a^{(k)}+k\geq k+1$, so that
after $k\leq a-1$ steps we get $b^{(k)}_{21}=b^{(k)}_{22}=0$ and
the procedure is completed. As a result we obtain the inequality
$$
\mu_{l,N}(B,d)\leq(a+1)d.
$$
Note that for $l=0$ (when all effective cycles are of the form
$d{\mathbb P}^N$) this estimate is precise: the equalities
$b_1=0$, $b_2=1$ mean that for a generic tuple $(f_1,\dots,f_N)\in
B$ the complete intersection
$$
\{f_1=\dots=f_{N-1}=0\}
$$
is a curve, non-singular at the point $o$. The condition of
tangency of order $j\leq N$ with that curve imposes at most $j$
independent conditions on the polynomial $f_N$. Therefore, the
equality
$$
\mu_{0,N}(a,d)=(a+1)d
$$
holds.\vspace{0.3cm}

%%%%%%%%%%%%%%%%%%%%%%%%%%%%%%%%%%%%%%%%%%%%%%%%%%%%%%%%%%%%%%%%%%%%%%%%%%%%%%%%%%%%
%%%%%%%%%%%%%%%%%%%%%%%%%%%%%%%%%%%%%%%%%%%%%%%%%%%%%%%%%%%%%%%%%%%%%%%%%%%%%%%%%%%%
%%%%%%%%%%%%%%%%%%%%%%%%%%%%%%%%%%% subsection 4.2

{\bf 4.2. Estimating the multiplicity in the case $b=2$.}
Similarly to Sec. 2.2, set $\mu_{l,N}(a,b;d)=m$, if for any
irreducible (bi-invariant) subvariety $B\subset{\cal
P}^{(n_1,n_2)}_N$ of codimension at most $a$ with
$\varepsilon(B)=b$ the inequality $\mu_{l,N}(B,d)\leq m$ holds,
and for at least one subvariety $B$ in that class the inequality
is an equality. The result of Sec. 4.1 can be represented as the
inequality
$$
\mu_{l,N}(a,1;d)\leq(a+1)d.
$$
Now let us obtain an upper bound for $\mu_{l,N}(a,2;d)$. Let $B$
be a subvariety with $\varepsilon(B)=2$. Applying the result of
Sec. 4.1, we get
$$
\mu_{l,N}(B,d)\leq(a-2)d+\mu_{l,N-1}(B_2,d_{21})+\mu_{l-1,N-1}(B_2,d_{22})
$$
for some $d_{21},d_{22}\in{\mathbb Z}_+$ with $d_{21}+d_{22}=d$.
For the parameters of the subvariety $B_2$ there are two options:
\begin{itemize}
\item either $\mathop{\rm codim}B_2\leq a-2$ and $\varepsilon(B_2)=2$,
\item or $\mathop{\rm codim}B_2\leq a-3$ and $\varepsilon(B_2)=1$.
\end{itemize}
In the second case we get the estimate
$$
\mu_{l,N}(B,d)\leq 2(a-2)d.
$$
In the first case one can go on with the process of reduction,
applying Proposition 3.3 to $B_2=B^{(1)}$. Arguing as in Sec. 4.1
(the computations are absolutely elementary and similar to those
performed in Sec. 4.1 and we do not give them here), we obtain the
following final result: for an even $a=2u$ the inequality
$$
\mu_{l,N}(a,2;d)\leq(2+u(u-1))d
$$
holds, for an odd $a=2u+1$ the inequality
$$
\mu_{l,N}(a,2;d)\leq(2+u^2)d
$$
holds.\vspace{0.1cm}

This procedure of obtaining explicit upper bounds for the numbers
$\mu_{l,N}(a,b;d)$ can be iterated, reducing the estimate for
$b=3$ to the already known formulas for $b=1,2$. However, as could
already see in the case $b=2$, the number of cases that require
separate consideration, starts to grow, and the formulas get
clumsier. Thus for small codimensions it is easier to obtain a
particular numerical value of the upper bound for
$\mu_{l,N}(a,b;d)$, whereas for higher values of $a$ we need a
less precise but manageable estimate.\vspace{0.3cm}

%%%%%%%%%%%%%%%%%%%%%%%%%%%%%%%%%%%%%%%%%%%%%%%%%%%%%%%%%%%%%%%%%%%%%%%%%%%%%%%%%%%%
%%%%%%%%%%%%%%%%%%%%%%%%%%%%%%%%%%%%%%%%%%%%%%%%%%%%%%%%%%%%%%%%%%%%%%%%%%%%%%%%%%%%
%%%%%%%%%%%%%%%%%%%%%%%%%%%%%%%%%%% subsection 4.3

{\bf 4.3. Small codimensions.} A simple observation that one can
make on the basis of the considerations for $b=1$ and $b=2$, is
the linearity of the obtained estimates in the degree $d$ and
their actual independence of the parameters $l,N$. Let
$U\subset{\mathbb Z}_+\times{\mathbb Z}_+$ be the set
$\{(a,b)\,|\,a\geq b^2\}$. Let us define by induction the function
$$
\overline{\mu}\colon U\to{\mathbb Z}_+,
$$
setting $\overline{\mu}(a,0)\equiv 1$, $\overline{\mu}(a,1)\equiv
a+1$, for $a<b(b+1)$
$$
\overline{\mu}(a,b)=2\overline{\mu}(a-(2b-1),b-1),
$$
for $a\geq b(b+1)$
$$
\overline{\mu}(a,b)=\overline{\mu}(a-(2b-1),b-1)+\mathop{\rm max}\{\overline{\mu}(a-(2b-1),b-1),
\overline{\mu}(a-b,b)\}.
$$
Propositions 3.1, 3.2 and 3.3 imply immediately\vspace{0.1cm}

{\bf Proposition 4.1.} {\it The following inequality holds:}
$$
\mu_{l,N}(a,b;d)\leq\overline{\mu}(a,b)d.
$$
\vspace{0.1cm}

For small values of $a$ the function $\overline{\mu}$ is easy to
compute by hand; it is also easy to write a computer program,
computing $\overline{\mu}$. Below we give the table of values
$\overline{\mu}(a,b)$ for $a\leq 36$, $b\leq 6$. The symbol $*$
means that the pait $(a,b)\not\in U$ and the value of the function
$\overline{\mu}$ is not defined. Already for those small values of
the codimension the speed of growth of the values
$\overline{\mu}(a,b)$ can be seen very well. In boldface we show
the maximum value $\overline{\mu}(a,b)$ for the given
$a$.\vspace{0.5cm}

\begin{tabular}{|l|c|c|c|c|c|c|c|c|c|c|c|c|c|c|c|}
\hline
$a$   & 1 & 2 & 3 & 4 & 5 & 6 & 7  & 8  & 9  & 10 & 11 & 12 & 13 & 14 & 15  \\
\hline
$b=0$ & 1 & 1 & 1 & 1 & 1 & 1 & 1  & 1  & 1  & 1  & 1  & 1  & 1  & 1  & 1   \\
\hline
$b=1$ & {\bf 2} & {\bf 3} & {\bf 4} & {\bf 5} & {\bf 6} & 7 & 8  & 9  & 10 &
11 & 12 & 13 & 14 & 15 & 16  \\
\hline
$b=2$ & * & * & * & 4 & 6 & {\bf 8} & {\bf 11} & {\bf 14} & {\bf 18} & {\bf 22} &
{\bf 27} & {\bf 32} & {\bf 38} & {\bf 44} & {\bf 51}  \\
\hline
$b=3$ & * & * & * & * & * & * & *  & *  & 8  & 12 & 16 & 22 & 28 & 36 & 44  \\
\hline
$b=4$ & * & * & * & * & * & * & *  & *  & *  & *  & *  & *  & *  & *  & *   \\
\hline
\end{tabular}
\vspace{0.5cm}

\begin{tabular}{|l|c|c|c|c|c|c|c|c|c|c|c|}
\hline
$a$   & 16 & 17 & 18 & 19 & 20  & 21  & 22  & 23  & 24  & 25  & 26  \\
\hline
$b=0$ & 1  & 1  & 1  & 1  & 1   & 1   & 1   & 1   & 1   & 1   & 1   \\
\hline
$b=1$ & 17 & 18 & 19 & 20 & 21  & 22  & 23  & 24  & 25  & 26  & 27  \\
\hline
$b=2$ & {\bf 58} & 66 & 74 & 83 & 92  & 102 & 112 & 123 & 134 & 146 & 158 \\
\hline
$b=3$ & 55 & {\bf 68} & {\bf 82} & {\bf 99} & {\bf 119} & {\bf 140} &
{\bf 165} & {\bf 193} & {\bf 223} & {\bf 257} & {\bf 295} \\
\hline
$b=4$ & 16 & 24 & 32 & 44 & 56  & 72  & 88  & 110 & 136 & 164 & 198 \\
\hline
$b=5$ & *  & *  & *  &  * &   * &   * &  *  & *   & *   & 32  & 48  \\
\hline
$b=6$ & *  & *  & *  & *  & *   & *   & *   & *   & *   & *   & *  \\
\hline
\end{tabular}
\vspace{0.5cm}

\begin{tabular}{|l|c|c|c|c|c|c|c|c|c|c|}
\hline
$a$   & 27  & 28  & 29  & 30  & 31  & 32  & 33  & 34  & 35  & 36  \\
\hline
$b=0$ & 1   & 1   & 1   & 1   & 1   & 1   & 1   & 1   & 1   & 1   \\
\hline
$b=1$ & 28  & 29  & 30  & 31  & 32  & 33  & 34  & 35  & 36  & 37  \\
\hline
$b=2$ & 171 & 184 & 198 & 212 & 227 & 242 & 258 & 274 & 291 & 308 \\
\hline
$b=3$ & {\bf 335} & {\bf 380} & {\bf 429} & {\bf 481} & {\bf 538} &
{\bf 600} & {\bf 665} & {\bf 736} & 812 & 892 \\
\hline
$b=4$ & 238 & 280 & 330 & 391 & 461 & 537 & 625 & 726 & {\bf 841} & {\bf 966} \\
\hline
$b=5$ & 64  & 88  & 112 & 144 & 176 & 220 & 272 & 328 & 396 & 476 \\
\hline
$b=6$ & *   & *   & *   & *   & *   & *   & *   & *   & *   & 64  \\
\hline
$b=7$ & *   & *   & *   & *   & *   & *   & *   & *   & *   & *   \\
\hline
\end{tabular}
\vspace{0.5cm}

Now let us consider the problem of obtaining a simple effective
upper bound for the multiplicities $\mu_{l,N}(B,d)$. From the
technical viewpoint, it is necessary to find a simple and visual
formalization of the procedure of estimating these numbers in
terms of the numbers $\mu_{l',N'}(B',d')$, where the varieties
$B'\subset{\cal P}^{(n'_1,n'_2)}_{N'}$ have a smaller value
$b'=b'_1+b'_2$, so that for the corresponding multiplicities the
upper bound can be assumed to be known. This gives an inductive
(in the parameter $b$) procedure of estimating the multiplicity,
realized below.\vspace{0.3cm}

%%%%%%%%%%%%%%%%%%%%%%%%%%%%%%%%%%%%%%%%%%%%%%%%%%%%%%%%%%%%%%%%%%%%%%%%%%%%%%%%%%%%
%%%%%%%%%%%%%%%%%%%%%%%%%%%%%%%%%%%%%%%%%%%%%%%%%%%%%%%%%%%%%%%%%%%%%%%%%%%%%%%%%%%%
%%%%%%%%%%%%%%%%%%%%%%%%%%%%%%%%%%% subsection 4.4

{\bf 4.4. A general method of estimating the multiplicity.} Let us
consider the four-letter alphabet $\{A,C_0,C_1,C_2\}$. We define a
procedure of constructing a certain set of words $W$ in this
alphabet and for each word $w\in W$ a certain irreducible
bi-invariant subvariety
$$
B[w]\subset{\cal P}^{(n_1(w),n_2(w))}_{N(w)}
$$
of the type $\tau(B[w])=((a_1(w),b_1(w)),(a_2(w),b_2(w))$ and full
codimension $a(w)=a_1(w)+a_2(w)$; as usual, set
$\varepsilon(B[w])=b(w)=b_1(w)+b_2(w)$. The length of the word $w$
we denote by the symbol $|w|\in{\mathbb Z}_+$. The length of the
empty word is equal to zero.\vspace{0.1cm}

Let $B\subset{\cal P}^{(n_1,n_2)}_N$ be an irreducible
bi-invariant variety of codimension $a$ and type $\tau(B)$. Set
$B[\emptyset]=B$. If $b(\emptyset)=0$, then we set
$W=\{\emptyset\}$: the procedure is thus complete. Assume that
$b(\emptyset)\geq 1$.\vspace{0.1cm}

The set of words $W$ and corresponding subvarieties $B[w]$, $w\in
W$, will be constructed in elementary steps: we will construct a
sequence of finite subsets $W_l$ of the set of all words,
$l=0,1,\dots$. The set $W_0=\{\emptyset\}$ is already constructed.
Assume that $W_0,\dots,W_l$ are constructed. If for every $w\in
W_l$ the equality $b(w)=0$ holds, then we set $W=W_l$, completing
the procedure. Otherwise, take any word $w\in W_l$ with $b(w)\geq
1$. If $b_2(w)\geq 1$, then we apply to the subvariety $B[w]$
(constructed at the previous step) part (i) of Proposition 3.3.
Set $w_1=wA$ and $w_2=wC_e$, where $e\in\{0,1,2\}$ is chosen in
the following way.\vspace{0.1cm}

Set $B[w_1]=(B[w])_1$ and $B[w_2]=(B[w])_2$ in the notations of
Proposition 3.3, (i). Now $e$ takes the value 0, 1 or 2, if the
subvariety $B[w_2]$ corresponds to the case (i), (ii) or (iii) of
Proposition 3.1, respectively. This determines the words $w_1$ and
$w_2$. The set $W_{l+1}$ is obtained from $W_l$ by removing the
word $w$ and adding the words $w_1,w_2$ of length $|w|+1$. The
irreducible bi-invariant subvarieties $B[w_i]$ were constructed
above, and this defines the values of all parameters ($n_i(w)$
etc.).\vspace{0.1cm}

If $b_2(w)=0$, then $b_1(w)\geq1$ and we apply to the subvariety
$B[w]$ part (ii) of Proposition 3.3. Now the words $w_1,w_2$ and
the subvarieties $B[w_i]$ are constructed word for word in the
same way as in the case $b_2(w)\geq 1$, replacing part (i) of
Proposition 3.3 by part (ii) of the same proposition and
Proposition 3.1 by Proposition 3.2.\vspace{0.1cm}

Now the procedure of constructing the set $W_{l+1}$ is determined
in a unique way. Obviously,
$$
\sharp W_l=l+1.
$$
Since when we replace a word $w\in W_l$ by the words $w_1,w_2$,
the codimension of the new subvarieties $B[w_i]$ gets strictly
smaller than the codimension of $B[w]$ (Propositions 3.1 and 3.2),
our procedure of constructing the sequence $\{W_l\}$ can not be
infinite. It is easy to see that for $\mathop{\rm codim}B=a$ we
get $\sharp W\leq 2^a$. Now Proposition 3.3 implies
immediately\vspace{0.1cm}

{\bf Proposition 4.2.} {\it The following inequality holds:}
$$
\mu_{l,N}(B,d)\leq d(\sharp W).
$$
\vspace{0.1cm}

{\bf Proof.} This follows from a more general fact:
\begin{equation}\label{oct11.12}
\mu_{l,N}(B,d)\leq
\sum_{w\in W_i}\sum^{{\rm min}\{l,|w|\}}_{j=0}\mu_{l-j,N-|w|}(B[w],d_j(w))
\end{equation}
for any $i=0,1,\dots$ and certain partitions
\begin{equation}\label{oct11.13}
d=\sum^{{\rm min}\{l,|w|\}}_{j=0}d_j(w)
\end{equation}
for every word $w\in W_i$. The inequality (\ref{oct11.12}) will be
shown by induction on the parameter $i=0,1,\dots$. For $i=0$ its
left hand side and right hand side are the same. Assume that
(\ref{oct11.12}) is shown for $i=0,\dots,e$. If $W_e=W_{e+1}$,
then $W_e=W$ and there is nothing more to prove. If $W_e\neq
W_{e+1}$, then the set of words $W_{e+1}$ is obtained from $W_e$
by removing some word $w\in W_e$ and ading two words $w_1$, $w_2$.
Therefore, in order to prove the inequality (\ref{oct11.12}) for
$i=e+1$ it is sufficient to show that
$$
\sum^{{\rm min}\{l,|w|\}}_{j=0}\mu_{l-j,N-|w|}(B[w],d_j[w])
$$
does not exceed the sum of similar expressions for $w_1,w_2$. This
is precisely what Proposition 3.3 claims, taking into account the
inequality (\ref{oct11.11}). This proves the estimate
(\ref{oct11.12}) for any $i$. Finally, if $W_i=W$, then for any
$w\in W$ we have $b(w)=0$, so that in the right hand side of the
inequality (\ref{oct11.12}) for each component we have the
equality
$$
\mu_{l-j,N-|w|}(B[w],d_j(w))=d_j(w),
$$
so that by the equality (\ref{oct11.13}) we get
$$
\mu_{l,N}(B,d)\leq d\sum_{w\in W}1,
$$
which is what we claimed. Q.E.D. for Proposition
4.2.\vspace{0.3cm}

%%%%%%%%%%%%%%%%%%%%%%%%%%%%%%%%%%%%%%%%%%%%%%%%%%%%%%%%%%%%%%%%%%%
%%%%%%%%%%%%%%%%%%%%%%%%%%%%%%%%%%%%%%%%%%%%%%%%%%%%%%%%%%%%%%%%%%%
%%%%%%%%%%%%%%%%%%%%%%%%%%%%%%%%%%% subsection 4.5

{\bf 4.5. Estimating the cardinality of the set of words.} We
write down the words in the following way:
$$
w=\tau_1\dots\tau_K,
$$
where $\tau_i\in\{A,C_0,C_1,C_2\}$. Now let
$$
\nu\colon\{A,C_0,C_1,C_2\}\to\{A,C\}
$$
be the map of the four-letter alphabet into the two-letter one,
given by the equalities $\nu(A)=A$, $\nu(C_i)=C$, and
$$
\nu\colon w=\tau_1\dots\tau_K\mapsto\bar{w}=\nu(\tau_1)\dots\nu(\tau_K)
$$
the corresponding map of the set of words. The following fact is
true.\vspace{0.1cm}

{\bf Lemma 4.1.} {\it For any $i=0,1,\dots$ the map $\nu|_{W_i}$
is injective. In particular, $\nu|_W$ is injective.}\vspace{0.1cm}

{\bf Proof.} A stronger fact is true: no word among the words
$\bar{w}=\nu(w)$, $w\in W_i$, is a left segment of any other word
in this set. (In particular, no two words are the same, which
means precisely the injectivity of the map $\nu|_{W_i}$.) The last
claim is easy to show by induction. The set $W_0$ consists of one
word, and for this set the claim is trivial. Assume that it is
shown for $W_i$, where $i=0,\dots,e$. If $W_{e+1}=W_e$, then there
is nothing to prove. If $W_{e+1}\neq W_e$, then $W_{e+1}$ is
obtained from $W_e$ by removing some word $w\in W_e$ and adding
two words $w_1=wA$ and $w_2=wC_{\alpha}$, where
$\alpha\in\{0,1,2\}$. For these words we have $\bar{w}_1=\bar{w}A$
and $\bar{w}_2=\bar{w}C$. Obviously, $\bar{w}_1$ and $\bar{w}_2$
are not left segments of each other and no word $\bar{w}'_1$ for
$w'\in W_e\backslash\{w\}$ is a left segment of $\bar{w}_1$ or
$\bar{w}_2$, as otherwise $\bar{w}'=\bar{w}_1$ or $\bar{w}_2$
(since $\bar{w}'$ is not a left segment of the word $\bar{w}$ by
the induction hypothesis), but then $\bar{w}$ would be a left
segment of the word $\bar{w}'$, contrary to the induction
hypothesis. In the trivial way $\bar{w}_1$ and $\bar{w}_2$ are no
left segments of any word $\bar{w}'$, since otherwise this would
have been true for $\bar{w}$ as well, contrary to the induction
hypothesis. Q.E.D. for the lemma.\vspace{0.1cm}

Thus we have reduced the problem of estimating the multiplicity
$\mu_{l,N}(B,d)$ to the problem of estimating the number of words
in the set $W$. As we pointed out above, $\sharp W\leq 2^a$, but
that estimate is too coarse for our purposes. We will control the
length of words $w\in W$ via the values of the parameters $a(w')$
and $b(w')=\varepsilon(B[w])$ for the left segments $w'$ of the
word $w$.\vspace{0.1cm}

{\bf Lemma 4.2.} (i) {\it If $\tau=A$ or $C_0$, then the
inequality
$$
a(w'\tau)\leq a(w')-(2b(w')-1)
$$
holds and $b(w'\tau)=b(w')-1$.}\vspace{0.1cm}

(ii) {\it If $\tau=C_1$ or $C_2$, then the inequality
$$
a(w'\tau)\leq a(w')-b(w')
$$
holds and $b(w'\tau)=b(w')$.}\vspace{0.1cm}

{\bf Proof.} This follows immediately from the inequalities of
Propositions 3.1 and 3.2, taking into account the obvious estimate
$b_i\leq n_i$, $i=1,2$. Q.E.D. for the lemma.\vspace{0.1cm}

Furthermore, the inequality (\ref{oct11.4}) implies the estimate
$$
a(w)\geq b_1(w)(n_2(w)+b_1(w))+b_2(w)b(w)\geq b^2(w)
$$
for every word $w$, in particular, for every word, which is a left
segment of any word in $W$.\vspace{0.1cm}

{\bf Example 4.1.} In terms of the formalism, developed above, let
us again consider the case $b=b(\emptyset)=1$. Here for any word
$w\in W_i$ we have the alternative: either $b(w)=0$ (and in that
case $w\in W$), or $b(w)=1$ (and in that case $a(w\tau)\leq
a(w)-1$ for any letter $\tau$), so that the set $W$ is of the form
$$
A,\,\, C_{i_1}A,\,\, C_{i_1}C_{i_2}A,\,\, \dots,\,\,
C_{i_1}C_{i_2}\dots C_{i_k}A,\,\, C_{i_1}
\dots C_{i_k}C_0,
$$
where $i_{\alpha}\in\{1,2\}$ and $k+1\leq a$. Therefore, $\sharp
W\leq a+1$, as we claimed above in Sec.~4.1.\vspace{0.1cm}

Let us come back to the general case.\vspace{0.1cm}

{\bf Proposition 4.3.} {\it The following inequality holds:}
$$
\sharp W\leq2^b\frac{a^b}{(b!)^2}.
$$
\vspace{0.1cm}

{\bf Proof.} For any word $w\in W$ by construction $b(w)=0$. Since
the letters $C_1,C_2$ do not change the value of the parameter
$b$, whereas the letters $A$ and $C_0$ decrease it by 1, we can
conclude that in the word $w$ there are precisely $b$ positions,
occupied by the letters $A$ and $C_0$. Let them be the positions
with the numbers
$$
m_1+1,\,\, m_1+m_2+2,\,\, \dots,\,\, m_1+m_2+\dots+m_b+b,
$$
$m_i\in{\mathbb Z}_+$. By Lemma 4.2, we get the inequality
$$
\begin{array}{cclcl}
0\leq a(w)\leq a & - & m_1b & - & (2b-1)-\\
& - & m_2(b-1) & - & (2(b-1)-1)-\\
&  & \dots & & \\
& - & m_i(b-(i-1)) & - & (2(b-(i-1))-1)-\\
&  & \dots & & \\
& - & m_b & - & 1 = \\
&  &  & = & a - b^2-\sum\limits^b_{i=1}m_i(b-(i-1)),
\end{array}
$$
so that $(m_1,\dots,m_b)$ is an arbitrary integer-valued point in
the polytope
$$
\Delta=\{x_1\geq 0,\dots,x_b\geq 0,bx_1+(b-1)x_2+\dots+x_b\leq
a-b^2\}\subset{\mathbb R}^b.
$$
Thus even if we assume that all possible distributions of the
letters $A$ and $C_0$ on the selected positions are realized by
the words $w\in W$ (in fact, that is not true: there are much
fewer words in $W$, see Remark 4.1), then the following inequality
holds:
$$
\sharp W\leq 2^b\cdot\sharp(\Delta\cap{\mathbb Z}^b).
$$
Now let us estimate the number of integer-valued points in
$\Delta$. In order to do that, consider a larger polytope
$$
\Delta^+=\{x_1\geq 0,\dots,x_b\geq 0,bx_1+\dots+
x_b\leq a-\frac{b(b-1)}{2}\}\subset{\mathbb R}^b.
$$
Obviously, $\Delta\subset\Delta^+$.\vspace{0.1cm}

{\bf Lemma 4.3.} {\it The following inequality holds:}
$$
\sharp(\Delta\cap{\mathbb Z}^b)\leq\mathop{\rm vol}(\Delta^+).
$$
\vspace{0.1cm}

{\bf Proof.} With each point $x=(x_1,\dots,x_b)\in{\mathbb R}^b$
we associate the unit cube
$$
\Gamma(x)=[x_1,x_1+1]\times[x_2,x_2+1]\times\dots
\times[x_b,x_b+1]\subset{\mathbb R}^b,
$$
the vertex of which with the least value of the sum of coordinates
$x_1+\dots+x_b$ is the point $x$. If $x\in\Delta$, then
$\Gamma(x)\subset\Delta^+$, since
$$
b+(b-1)+\dots+1+a-b^2=a-\frac{b(b-1)}{2}.
$$
Therefore,
$$
\sharp(\Delta\cap{\mathbb Z}^b)=
\sum_{x\in\Delta\cap{\mathbb Z}^b}\mathop{\rm vol}(\Gamma(x))=
\mathop{\rm vol}\left(\bigcup_{x\in\Delta\cap{\mathbb Z}^b}\Gamma(x)\right)
\leq\mathop{\rm vol}(\Delta^+),
$$
as we claimed. Q.E.D. for the lemma.\vspace{0.1cm}

Computing the volume of the polytope $\Delta^+$ and applying the
Stirling formula, we get the estimate
\begin{equation}\label{9april12.10}
\sharp W\leq 2^b\frac{(a-\frac{b(b-1)}{2})^b}{(b!)^2}=
\frac{1}{2\pi be^{\theta/6b}}\left(\frac{2a-b(b-1)}{b^2}e^2\right)^b
\end{equation}
for some $0<\theta<1$ (here $e$ is the base of the natural
logarithm), so that the more so,
\begin{equation}\label{23april12.1}
\sharp W\leq u_b= \frac{1}{2\pi
b}\left(\frac{2a-b(b-1)}{b^2}e^2\right)^b.
\end{equation}
Recall that $b\in\{1,\dots,[\sqrt{a}]\}$. To obtain an effective
bound for the number $\sharp W$ let us study the behaviour of the
sequence $u_b$ for those values of $b$.\vspace{0.1cm}

{\bf Lemma 4.4.} {\it The sequence $u_b$ is increasing, provided
that the following inequality holds:}
\begin{equation}\label{9april12.11}
2a-b(b-1)\geq \frac52b^2.
\end{equation}
\vspace{0.1cm}

{\bf Proof.} Write down
\begin{equation}\label{9april12.12}
\frac{u_{b+1}}{u_b}=\frac{1}{1+\frac{1}{b}}
\frac{e^2}{\left(1+\frac{1}{b}\right)^{2b}}
\frac{1}{\left(1+\frac{2b}{2a-b(b+1)}\right)^b}\frac{2a-b(b+1)}{(b+1)^2}.
\end{equation}
Assume first that $b\geq 9$. If the numbers $a$ and $b$ satisfy
the inequality $2a-b(b+1)\geq \frac52 (b+1)^2$ (that is, the
inequality (\ref{9april12.11}) for $b+1$), then the denominator of
the third factor in the right hand side can be estimated from
above as follows:
$$
\left(1+\frac{2b}{2a-b(b+1)}\right)^b\leq
\left(1+\frac{4}{5}\frac{1}{b}\right)^b<e^{\frac45}.
$$
The second factor in the right hand side of the inequality
(\ref{9april12.12}) is strictly higher than one, and the fourth is
not smaller than $\frac52$. As a result we get:
$$
\frac{u_{b+1}}{u_b}>\frac{9}{10}\cdot \frac52\cdot e^{-\frac45}>1,
$$
which is what we need. For the smaller values $b\leq 8$ the second
and third factors in the right hand side of the inequality
(\ref{9april12.12}) could be estimated more precisely, and
elementary computations with some use of a computer complete the
proof of the lemma.\vspace{0.1cm}

{\bf Corollary 4.1.} {\it For $a\geq 17$ the value $b_{\rm
max}\in\{1,\dots,[\sqrt{a}]\}$, on which the maximum of the
sequence $u_b$ is attained, satisfies the inequality}
$$
2a-b_{\rm max}(b_{\rm max}-1)\leq \frac53 a.
$$
\vspace{0.1cm}

{\bf Proof.} By the previous lemma, the value $b_{\rm max}$
satisfies the inequality
$$
2a-b_{\rm max}(b_{\rm max}+1)\leq \frac52 b_{\rm max}^2
$$
(otherwise, the next element of the sequence $u_b$ would be
higher). Now elementary computations complete the proof of the
corollary.\vspace{0.1cm}

{\bf Corollary 4.2.} (i) {\it For $a\geq 17$ the following
estimate holds:}
$$
\sharp W\leq v_b=\frac{1}{2\pi b}\left(\frac{5a}{3b^2}e^2\right)^b.
$$
\vspace{0.1cm}

(ii) {\it For any $a$ the following estimate holds:}
$$
\sharp W\leq w_b=\frac{1}{2\pi b}\left(\frac{2a}{b^2}e^2\right)^b.
$$
\vspace{0.1cm}

{\bf Proof.} The claim (ii) follows immediately from the
inequality (\ref{23april12.1}), the claim (i) from the inequality
(\ref{23april12.1}), taking into account the previous
corollary.\vspace{0.1cm}

{\bf Corollary 4.3.} (i) {\it For $a\geq 17$ the following
estimate holds:}
$$
\mu_{l,N}(a,d)\leq
\frac{e^2}{2\pi[\sqrt{a}]}\left(\frac53 e^2\right)^{[\sqrt{a}]}.
$$
\vspace{0.1cm}

(ii) {\it For any $a$ the following estimate holds:}
$$
\mu_{l,N}(a,d)\leq
\frac{e^2}{2\pi[\sqrt{a}]}\left(2e^2\right)^{[\sqrt{a}]}.
$$
\vspace{0.1cm}

{\bf Proof.} The arguments are identical in both cases, the only
difference is which of the two claims of Corollary 4.2 is
used.\vspace{0.1cm}

Let us prove part (i). Arguing as in the proof of Lemma 4.4, we
conclude that the sequence $v_b$ is increasing. Therefore, its
maximum is attained at $b=[\sqrt{a}]$. Since
$$
a<(b+1)^2=b^2+2b+1,
$$
we get the inequality
$$
\left(\frac{a}{b^2}\right)^b\leq \left(1+\frac{2}{b}\right)^b<e^2,
$$
whence immediately follows the claim (i). The second part is shown
in word for word the same way. Q.E.D.\vspace{0.1cm}

It is easy to see that the claims of Theorems 4 and 5 are
contained in the claims of Corollaries 4.2 and 4.3, taking into
account the formula (\ref{23april12.1}).\vspace{0.1cm}

Q.E.D. for Theorems 4 and 5.\vspace{0.1cm}

{\bf Remark 4.1.} As we see from the given proof, the estimate
obtained above is not optimal and can be essentially improved. For
$b\approx\sqrt{a}$ we have $2a-b(b-1)\approx a$, so that in the
inequality of Corollary 4.3, (ii), the expression $(2e^2)$ can be
replaced by $e^2$. Furthermore, when proving Proposition 4.3, we
took into account all possible tuples of positions
$(m_1,\dots,m_b)$ and all possible ways of putting the letters $A$
and $C_0$ on the $b$ positions. However, since in the set of words
$\overline{W}=\nu(W)$ in the two-letter alphabet $\{A,C\}$ no word
is a left segment of another word and the map $\nu\colon
W\to\overline{W}$ is one-to-one, for a fixed way of putting the
letters $A$ and $C_0$ on $b$ positions, when at least two letters
$C_0$ are neighbours, not all tuples
$(m_1,\dots,m_b)\in\Delta\cap{\mathbb Z}^b$ are realized, as two
distinct words $w_1\neq w_2$, $\{w_1,w_2\}\subset W$, can not
differ only on a segment that consists of the letters
$C_0,C_1,C_2$. The problem of getting the precise upper bound for
the numbers $\mu_{l,N}(a,d)$, at least in the asymptotic sense,
remains open.\vspace{0.5cm}

%%%%%%%%%%%%%%%%%%%%%%%%%%%%%%%%%%%%%%%%%%%%%%%%%%%%%%%%%%%%%%%%%%%%%%%%%%%%%%%%%%%%
%%%%%%%%%%%%%%%%%%%%%%%%%%%%%%%%%%%%%%%%%%%%%%%%%%%%%%%%%%%%%%%%%%%%%%%%%%%%%%%%%%%%
%%%%%%%%%%%%%%%%%%%%%%%%%%%%%%%%%%%%%%%%%%%%%%%%%%%%%%%%%%%%%%%%%%%%%%%%%%%%%%%%%%%%
%%%%%%%%%%%%%%%%%%%%%%%%%%%%%%%%%%%%%%%%%%%%%%%%%%%%%%%%%%%%%%%%%%%%%%%%%%%%%%%%%%%%
%%%%%%%%%%%%%%%%%%%%%%%%%%%%%%%%%%% SECTION 5

\begin{center}
\bf \S 5. Global multiplicities. Proof of Theorem 2.
\end{center}
\vspace{0.3cm}

In this section, using the theory developed in \S\S 2-4, we prove
Theorem 2. Taking into account Theorem 2 in the paper
\cite{Pukh11}, two facts require a proof: the linear independence
of the directions the lines passing through the point $o\in V$
(the last requirement in the condition (R1)), and that the
condition (R3) is satisfied at every point $o\in V$ on a Zarisky
generic complete intersection $V$ of the type $2^{k_1}3^{k_2}$. It
is not hard to prove the linear independence: it is sufficient to
estimate the codimension of the sets of tuples of polynomials,
which either have a positive-dimensional set of solutions or a
finite set set of linearly dependent solutions. This is done in
Sec. 5.1.\vspace{0.1cm}

In Sec. 5.2-5.4 we globalize the constructions and results of the
local theory, developed in \S\S 2-4: define the global
multiplicities, reduce the problem of their estimating to the
similar problem for the local multiplicities and, finally, obtain
the necessary estimates for the global
multiplicities.\vspace{0.1cm}

In Sec. 5.5 we complete the proof of Theorem 2.\vspace{0.3cm}

{\bf 5.1. Tuples of polynomials with a positive-dimensional set of
solutions.} As always, the symbol ${\cal P}_{i,N}$ denotes the
space of homogeneous polynomials of degree $i$ in $N$ variables,
we identify ${\cal P}_{i,N+1}$ and $H^0({\mathbb P}^N,{\cal
O}_{{\mathbb P}^N}(i))$. Let
$$
{\cal H}^{(n_1,n_2)}_N=\prod^{n_1}_{i=1}{\cal
P}_{2,N}\times\prod^{n_2}_{i=1}{\cal P}_{3,N}
$$
be the space of all tuples
$(f_1,\dots,f_{n_1},f_{n_1+1},\dots,f_{n_1+n_2})$, where the first
$n_1$ polynomials are quadratic, the next $n_2$ ones are cubic.
(This is the global analog of the space ${\cal
P}^{(n_1,n_2)}_{N}$, introduced in \S 2.) In this section we
consider the space ${\cal H}^{(n_1,n_2)}_{N+1}$ for $n_1+n_2=N$
and $N+1$. The symbol
$$
Z(f_1,\dots,f_{n_1+n_2})=Z(f_*)
$$
denotes the closed subscheme, defined by the tuple of polynomials
$f_1,\dots, f_{n_1+n_2}$, and the symbol $|Z(f_*)|$ stands for the
closed set $\{f_1=\dots=f_{n_1+n_2}=0\}\subset{\mathbb P}^N$. Let
$$
Y_{\infty}= \{(f_1,\dots,f_{n_1+n_2})\,|\,\mathop{\rm
dim}Z(f_*)\geq 1\}\subset {\cal H}^{(n_1,n_2)}_{N+1}
$$
be the closed subset of tuples, the zeros of which have an
``incorrect'' dimension (for a generic tuple for $n_1+n_2=N$ the
set $Z(f_*)$ is zero-dimensional, for $n_1+n_2=N+1$ it is empty).
Furthermore, define $Y_{\rm line}\subset Y_{\infty}$ as the set,
consisting of such tuples $(f_*)$, that the set $Z(f_*)$ contains
a line in ${\mathbb P}^N$. It is easy to compute that for
$n_1+n_2=N$
$$
\mathop{\rm codim}Y_{\rm line}=n_1+2n_2+2,
$$
and for $n_1+n_2=N+1$
$$
\mathop{\rm codim}Y_{\rm line}=n_1+2n_2+4
$$
(the codimension is in both cases with respect to the space ${\cal
H}^{(n_1,n_2)}_{N+1}$). Set
$Y'_{\infty}=\overline{Y_{\infty}\backslash Y_{\rm line}}$ to be
the union of all irreducible components of the set $Y_{\infty}$,
different from $Y_{\rm line}$. Now we have\vspace{0.1cm}

{\bf Proposition 5.1.} {\it The irreducible closed set $Y_{\rm
line}$ is a component of maximal dimension of the closed set
$Y_{\infty}$}:
$$
\mathop{\rm dim}Y_{\rm line}\geq\mathop{\rm dim}Y'_{\infty}.
$$
\vspace{0.1cm}

{\bf Proof.} The irreducibility of the set $Y_{\rm line}$ is
obvious. The claim of the proposition will be proved by the method
developed in \cite{Pukh01}. Consider first the case $n_1+n_2=N$.
Let $Y_i\subset {\cal H}^{(n_1,n_2)}_{N+1}$ be the set of such
tuples $(f_*)$, that\vspace{0.1cm}

(1) $\mathop{\rm codim}\{f_1=\dots=f_i=0\}=i$,\vspace{0.1cm}

(2) there is a component $B$ of the set $Z(f_1,\dots,f_i)$, which
for $i=N-1$ is not a line in ${\mathbb P}^N$ and on which the
polynomial $f_{i+1}$ vanishes identically. It is sufficient to
show that
$$
\mathop{\rm codim}\overline{Y_i}>n_1+2n_2+2
$$
for $i=1,\dots,N-1$.\vspace{0.1cm}

Following \cite[\S3]{Pukh01}, represent $Y_i$ as the union
$$
Y_i=Y_{i,0}\cup Y_{i,1}\cup\dots\cup Y_{i,i}
$$
of smaller subsets $Y_{i,e}$, $e\in\{0,\dots,i\}$, defined by the
condition: the set $B$ has codimension $e\in{\mathbb Z}_+$ in its
linear span $<B>$. Let us consider, first of all, the case $e=0$,
that is, $B\subset{\mathbb P}^N$ is a linear subspace, $i\leq
N-2$. Since $f_j|_B\equiv 0$ for $j=1,\dots,i+1$ and
$d_j=\mathop{\rm deg}f_j\geq 2$, we get the estimate
$$
\mathop{\rm codim}Y_{i,0}\geq(i+1)\frac{(N-i+1)(N-i+2)}{2}-i(N-i+1)=
$$
$$
=\frac{N-i+1}{2}(i(N-i)+N-i+2)\geq 3N,
$$
which is what we need (the minimum of the right hand side is
attained at $i=N-2$, which corresponds to the set of such tuples
$(f_*)$ that $Z(f_*)$ contains a plane).\vspace{0.1cm}

Therefore, we may assume that $e\geq 1$. In that case consider the
restrictions $f_j|_{<B>}$, $j=1,\dots,i$. Recall the following
\cite{Pukh01} \vspace{0.1cm}

{\bf Definition 5.1.} Let $h_1,\dots,h_m$ be homogeneous
polynomials of degree $\geq 2$ on the projective space $\Pi$ of
dimension $\mathop{\rm dim}\Pi\geq m+1$. An irreducible subvariety
$C\subset\Pi$ such that $<C>=\Pi$ and $\mathop{\rm codim}C=m$ is
called an {\it associated subvariety of the sequence} $(h_*)$, if
there exists a chain of irreducible subvarieties $C_j\subset\Pi$,
$j=0,\dots,m$, satisfying the following properties:
\begin{itemize}

\item $C_0=\Pi$,

\item for each $j=0,\dots,m-1$ the subvariety $C_{j+1}$ is an
irreducible component of the closed algebraic set
$\{h_{j+1}=0\}\cap C_j$, and moreover, $h_{j+1}|_{C_j}\not\equiv
0$, so that $\mathop{\rm codim}_{\Pi}C_j=j$ for all $j$,

\item $C_m=C$.
\end{itemize}

If a sequence $(h_*)$ has an associated subvariety, it is said to
be {\it good}.\vspace{0.1cm}

Furthermore, the following claim is true.\vspace{0.1cm}

{\bf Lemma 5.1.} (i) {\it The property of being good is an open
property in the space of sequences.}\vspace{0.1cm}

(ii) {\it A good sequence $(h_*)$ can have, at most,
$$
\left[\frac{1}{m+1}\prod^m_{j=1}\mathop{\rm deg}h_j\right]
$$
associated subvarieties.}\vspace{0.1cm}

{\bf Proof.} This is Lemma 4 in \cite{Pukh01}.\vspace{0.1cm}

As it was shown in \cite[p.73]{Pukh01}, one can find among the
polynomials $f_j$, $j=1,\dots,i$, such $e$ polynomials
$f_{j_1},\dots,f_{j_e}$ that the sequence
$(f_{j_1}|_{<B>},\dots,f_{j_e}|_{<B>})$ is good and $B$ is one of
its associated subvarieties. Besides, as it was shown in
\cite[p.72]{Pukh01}, for a fixed irreducible subvariety $C$ in the
projective space $\Pi$, such that $<C>=\Pi$, the requirement
$h|_C\equiv 0$ imposes on the polynomial $h$ at least
$$
\mathop{\rm deg}h\cdot\mathop{\rm dim}\Pi+1
$$
independent conditions. For that reason, fixing the subspace $<B>$
and the polynomials $f_{j_1},\dots,f_{j_e}$, we obtain that the
requirement $f_j|_B\equiv 0$ imposes on the polynomials $f_j$,
where
$$
j\in\{i+1\}\cup(\{1,\dots,i\}\backslash\{j_1,\dots,j_e\})
$$
at least
$$
(i-e+1)(2(N-i+e)+1)
$$
independent conditions. Subtracting the dimension of the
Grassmanian of $(N-i+e)$-subspaces in ${\mathbb P}^N$, we get
$$
\mathop{\rm codim}Y_{i,e}\geq(i-e+2)(N-i+e)+1.
$$
It is easy to check that the minimum of the right hand side is
$2N+1$. The inequality
$$
2N+1\geq n_1+2n_2+2
$$
is always true, except for the only case $n_1=0$, $n_2=N$.
However, estimating the number of independent conditions above, we
assumed that $\mathop{\rm deg}f_j=2$. If all the polynomials are
cubic ones, the estimate becomes essentially stronger:
$$
\mathop{\rm codim}Y_{i,e}\geq(2i-2e+3)(N-i+e)+1\geq 3N+1,
$$
which completes the proof of the proposition in the case
$n_1+n_2=N$.\vspace{0.1cm}

Consider the case $n_1+n_2=N+1$. Let $(f_1,\dots,f_{N+1})$ be a
tuple of general position in some irreducible component $Q$ od the
set $Y'_{\infty}$. Then there are two options: either $\mathop{\rm
dim}Z(f_1,\dots,f_N)=1$ and $f_{N+1}$ vanishes on some irreducible
component $C$ of the set $Z(f_1,\dots,f_N)$, and moreover the
curve $C$ is not a line, or $\mathop{\rm dim}Z(f_1,\dots,f_N)\geq
2$ and $f_{N+1}$ is an arbitrary polynomial.\vspace{0.1cm}

Assume that the first option takes place. According to what was
proved above, the codimension of the set
$$
\{(f_1,\dots,f_N)\,|\,(f_1,\dots,{f_{N+1}})\in Q\}\subset{\cal
H}^{(n_1,n_2-1)}_{N+1}
$$
is not smaller than $n_1+2(n_2-1)+2$, if $n_1\leq N$ (that is, if
at least one polynomial is a cubic one), and $n_1+1$, if
$n_1=N+1$. It is easy to see that the condition $f_{N+1}|_C\equiv
0$, where $C$ is a curve, the linear span of which is of dimension
at least 2, imposes on the polynomial $f_{N+1}$ at least 5
condition. Therefore, $\mathop{\rm codim} Q
>\mathop{\rm codim}Y_{\rm line}$, as we claimed.\vspace{0.1cm}

Assume that the second option takes place. Here we may assume that
$\mathop{\rm dim}Z(f_1,\dots,f_N)=2$. The polynomial $f_{N+1}$ is
an arbitrary one, so swapping $f_{N+1}$ and some polynomial $f_i$,
$i\leq N$, we get the situation considered above. Q.E.D. for
proposition 5.1.\vspace{0.1cm}

Now consider the set
$$
Y^{\Delta}_i\subset{\cal H}^{(n_1,n_2)}_{N+1}\backslash
Y_{\infty},
$$
consisting of such tuples $(f_1,\dots,f_{N+1})$, that there are
$i+1$ distinct points $p_1,\dots,p_{i+1}\in Z(f_*)$, such that
$$
\mathop{\rm dim}<p_1,\dots,p_{i+1}>\leq i-1.
$$
\vspace{0.1cm}

The following claim is true.\vspace{0.1cm}

{\bf Proposition 5.2.} {\it The equality}
$$
\mathop{\rm codim}Y^{\Delta}_i=N+2
$$
{\it holds.}\vspace{0.1cm}

{\bf Proof.} For convenience set $Y^{\Delta}_1=\emptyset$. It is
sufficient to show that $\mathop{\rm codim}(Y^{\Delta}_i\backslash
Y^{\Delta}_{i-1})=N+2$, that is, we may assume that among the
points $p_*$ any $i$ points are linearly independent, so that
$\Lambda=<p_*>$ is a subspace of dimension $i-1$. Let us fix
$\Lambda$. It is easy to see that the linear conditions
$f_j(p_i)=0$ are linearly independent, so that the set
$$
\{(f_1,\dots,f_{N+1})\,|\,\{p_1,\dots,p_{i+1}\}\subset |Z(f_*)|\}
$$
has codimension $(i+1)(N+1)$. The set of tuples
$\{p_1,\dots,p_{i+1}\}\subset\Lambda$ has dimension $(i+1)(i-1)$,
whereas the dimension of the projective Grassmanian is $i(N-i+1)$.
As a result we get
$$
\mathop{\rm codim}(Y^{\Delta}_i\backslash Y^{\Delta}_{i-1})=(i+1)(N-i+2)-i(N-i+1)=N+2,
$$
as we claimed. Q.E.D. for the proposition.\vspace{0.1cm}

Propositions 5.1 and 5.2 immediately imply\vspace{0.1cm}

{\bf Proposition 5.3.} {\it For $a\leq N+1$  the equality}
$$
\mathop{\rm codim}\{(f_1,\dots,f_{N+1})\,|\,\sharp |Z(f_*)|\geq
a\}=a
$$
{\it holds.}\vspace{0.1cm}

{\bf Proof.} One may assume that the set $|Z(f_*)|$ is finite and
consists of $a$ linearly independent points. In that case the
claim of the proposition is obvious. Q.E.D.\vspace{0.3cm}

%%%%%%%%%%%%%%%%%%%%%%%%%%%%%%%%%%%%%%%%%%%%%%%%%%%%%%%%%%%%%%%%%%
%%%%%%%%%%%%%%%%%%%%%%%%%%%%%%%%%%%%%%%%%%%%%%%%%%%%%%%%%%%%%%%%%%
%%%%%%%%%%%%%%%%%%%%%%%%%%%%%%%%%%% subsection 5.2

{\bf 5.2. The global multiplicities: the definitions and setting
up the problem.} We start to globalize the local theory, developed
in \S\S 2-4. Now, in order to distinguish between the local and
global multiplicities, for the local multiplicities of the type
$\mu_{l,N}(a,b;d)$ we write $\mu_{l,N}^{\rm local}(a,b;d)$ etc.
Consider the space
$$
{\cal H}^{(n_1,n_2)}_{N+1}= ({\cal P}_{2,N+1})^{n_1}\times({\cal
P}_{3,N+1})^{\times n_2}
$$
of tuples $(f_1,\dots,f_{n_1},f_{n_1+1},\dots,f_{n_1+n_2})$,
consisting of $n_1$ quadratic and $n_2$ cubic polynomials, which
we see as polynomials on ${\mathbb P}^N$. Assume that
$n_1+n_2=N+1$. Let
$$
\Sigma_1(f_*)=<f_1,\dots,f_{n_1}>
$$
be the linear system, generated by the quadratic polynomials, and
$$
\Sigma_2(f_*)=<f_{n_1+1},\dots,f_{n_1+n_2}>+\,\Sigma_1{\cal
P}_{1,N+1}
$$
be the linear system of cubic polynomials, generated by all
polynomials $f_*$.\vspace{0.1cm}

By the symbol $Y_{\infty}$ we denote the closed subset of the
space ${\cal H}^{(n_1,n_2)}_{N+1}$, consisting of the tuples
$(f_*)$ with the zero set of positive dimension: $\mathop{\rm
dim}Z(f_*)\geq 1$. According to Proposition 5.1, we have
$\mathop{\rm codim}Y_{\infty}=n_1+2n_2+4$. Let
$$
(f_1,\dots,f_{N+1})\in{\cal H}^{(n_1,n_2)}_{N+1}\backslash
Y_{\infty}
$$
be an arbitrary tuple.\vspace{0.1cm}

We write $\Sigma_1=\Sigma_1(f_*)$ and $\Sigma_2=\Sigma_2(f_*)$, if
it is clear, which tuple is meant. Consider the set
$$
\Sigma^{n_1}_1\times\Sigma^{N-n_1}_2\subset
({\cal P}_{2,N+1})^{\times n_1}\times
({\cal P}_{3,N+1})^{\times n_2}.
$$
Let
$(f^{\sharp}_1,\dots,f^{\sharp}_N)\in\Sigma^{n_1}_1\times\Sigma^{N-n_1}_2$
be a tuple of general position. For an arbitrary effective cycle
$R$ of pure codimension $l\in{\mathbb Z}_+$ on ${\mathbb P}^N$ we
define the {\it global effective multiplicity} by the equality
$$
\mu_{\rm total}((f_*);R)=
\sum_{p\in|Z(f_*)|}\mu_{\rm local}((f_*);R,p)
$$
(recall that the set $|Z(f_*)|$ is finite, since the tuple $(f_*)$
has been chosen outside the subset $Y_{\infty}$), where the local
multiplicity at the point $p$ is meant in the sense of Sec. 2.1:
$$
\mu_{\rm local}((f_*);R,p)= \mathop{\rm
mult}\nolimits_p(\{f^{\sharp}_{l+1}=0\}\circ
\dots\circ\{f^{\sharp}_N=0\}\circ R).
$$
(Obviously, in the right hand side in brackets we have a
zero-dimensional cycle, so that $\mathop{\rm mult}_p$ is the
multiplicity of the point $p$ in that cycle.) If $R\subset{\mathbb
P}^N$ is an irreducible subvariety of codimension $l$, then
$$
\mu_{\rm local}((f_*);R,p)=\mathop{\rm dim}{\cal
O}_{p,R}/(f^{\sharp}_{l+1},\dots,f^{\sharp}_N).
$$
\vspace{0.1cm}

Furthermore, let $Y^{\Delta}\subset{\cal H
}^{(n_1,n_2)}_{N+1}\backslash Y_{\infty}$ be the set of such
tuples $(f_*)$, that in the finite set of points $|Z(f_*)|$ one
can choose a linearly dependent subset. By Proposition 5.2,
$$
\mathop{\rm codim}\overline{Y^{\Delta}}=N+2.
$$
Therefore, we get the presentation
$$
{\cal H
}^{(n_1,n_2)}_{N+1}\backslash(Y_{\infty}\cup\overline{Y}^{\Delta})=
Y_0\sqcup Y_1\sqcup \dots\sqcup Y_{N+1},
$$
where $Y_a$ is a constructive set, such that for every tuple
$(f_*)\in Y_a$ the set of zeros $|Z(f_*)|$ consists of precisely
$a\in\{0,\dots,N+1\}$ linearly independent points. The closures
$\overline{Y_a}$ are irreducible, the sets $Y_a$ are open in their
closures and, by Proposition 5.3, $\mathop{\rm
codim}\overline{Y_a}=a$. In particular, $Y_0\subset{\cal
H}^{(n_1,n_2)}_{N+1}$ is an open subset.\vspace{0.1cm}

As in \S2, the symbol ${\mathbb H}_{l,N}(d)$ stands for the Chow
variety, parametrizing effective cycles of pure codimension
$l\in{\mathbb Z}_+$ and degree $d\geq 1$ on the projective space
${\mathbb P}^N$.\vspace{0.1cm}

Define the subset
$$
{\cal Y}^{\rm total}_{l,N,i}(m,d)\subset
Y_i\times{\mathbb H}_{l,N}(d)
$$
by the condition $\mu_{\rm total}((f_*);R)\geq m\in{\mathbb Z}_+$.
This is a closed algebraic subset (for a given value
$i\in\{0,\dots,N+1\}$). By one and the same symbol $\pi_{\cal P}$
we denote the projection of the direct product $Y_i\times{\mathbb
H}_{l,N}(d)$ onto $Y_i$. Since the Chow varieties are projective,
the image
$$
Y^{\rm total}_{l,N,i}(m,d)=
\pi_{\cal P}({\cal Y}^{\rm total}_{l,N,i}(m,d))\subset Y_i
$$
is a closed subset. It is clear that $Y_i$ and $Y^{\rm
total}_{l,N,i}(m,d)$ are all invariant with respect to the action
of the group $GL_{N+1}({\mathbb C})$ of linear changes of
coordinates. Let us consider the problem of estimating the
codimension
$$
\mathop{\rm codim}(Y^{\rm total}_{l,N,i}(m,d)\subset Y_i).
$$
\vspace{0.1cm}

Apart from the group $GL_{N+1}({\mathbb C})$, on the space ${\cal
H}^{(n_1,n_2)}_{N+1}$ acts the group $G^*(n_1,n_2)$. This action
is similar to the action of the group $G(n_1,n_2)$ in the local
situation. More precisely, the elements $g\in G^*(n_1,n_2)$ are
triples
$$
g=(A_{11}\in GL_{n_1}({\mathbb C}),A_{22}\in GL_{n_2}({\mathbb
C}), A_{12}\in \mathop{\rm Mat}\nolimits_{(n_1,n_2)}({\cal
P}_{1,N+1})),
$$
here $g(f_1,\dots,f_{N+1})=(f^g_1,\dots,f^g_{N+1})$, where
$$
(f^g_1,\dots,f^g_{n_1})=(f_1,\dots,f_{n_1})A_{11}
$$
and
$$
(f^g_{n_1+1},\dots,f^g_{n_1+n_2})=
(f_{n_1+1},\dots,f_{n_1+n_2})A_{22}+(f_1,\dots,f_{n_1})A_{12}.
$$
\vspace{0.1cm}

Subsets that are invariant both with respect to linear changes of
coordinates and the action of the group $G^*(n_1,n_2)$, are, as in
the local case, said to be {\it bi-invariant}.\vspace{0.1cm}

It is obvious that all the sets $Y_{\infty},Y^{\Delta},Y_i$ and
$Y^{\rm total}_{l,N,i}(m,d)$ are bi-invariant.\vspace{0.1cm}

By construction, the set $Y^{\rm total}_{l,N,i}(m,d)$ consists of
such tuples of polynomials $(f_1,\dots,f_{N+1})\in Y_i$, that
there exists an effective cycle $R$ of pure codimension $l$ and
degree $d$ on ${\mathbb P}^N$, satisfying the inequality $\mu_{\rm
total}((f_*);R)\geq m$. Similar to the local case, we re-formulate
the problem of estimating the codimension of the set $Y^{\rm
total}_{l,N,i}(m,d)$: we will maximize the multiplicity $m$ for a
fixed codimension. More precisely, let $B\subset Y_i$ be a
bi-invariant irreducible subvariety. Set
$$
\mu_{\rm total}(B,d)= \mathop{\rm max}\{m\,|\,B\subset Y^{\rm
total}_{l,N,i}(m,d)\}
$$
(the indices $l,N$ are omitted in the left hand side to simplify
the notations). Explicitly,
$$
\mu_{\rm total}(B,d)= \mathop{\rm max}\limits_R\{\mu_{\rm
total}((f_*);R)\},
$$
the maximum is taken over all effective cycles of degree $d$ and
codimension $l$, and, as usual, $(f_*)\in B$ is a generic tuple.
Finally, set $\mu_{\rm total}(i,a,d)=m$, if there exists an
irreducible bi-invariant subvariety $B\subset Y_i$ of codimension
at most $a$, such that $\mu_{\rm total}(B,d)=m$, and such a
variety does not exist for $m+j$, $j\geq 1$. Assuming the
codimension $a$ to be fixed, let us estimate the multiplicity $m$
from above.\vspace{0.3cm}

%%%%%%%%%%%%%%%%%%%%%%%%%%%%%%%%%%%%%%%%%%%%%%%%%%%%%%%%%%%%%%%%%
%%%%%%%%%%%%%%%%%%%%%%%%%%%%%%%%%%%%%%%%%%%%%%%%%%%%%%%%%%%%%%%%%
%%%%%%%%%%%%%%%%%%%%%%%%%%%%%%%%%%% subsection 5.3

{\bf 5.3. The local and global type of a subvariety of tuples of
polynomials.} Fix an irreducible bi-invariant subvariety $B\subset
Y_r$, $r\geq 1$. For a generic (and arbitrary) tuple $(f_*)\in B$
the set theoretic intersection $|Z(f_*)|$ consists of precisely
$r$ linearly independent points $p_1,\dots,p_r\in{\mathbb P}^N$
(depending on the tuple $(f_*)$, of course). For a point $p\in
|Z(f_*)|$ set $\varepsilon_p(f_*)=b$, where
$$
\mathop{\rm rk}(df_1(p),\dots,df_{N+1}(p))=N-b.
$$
Set $b_i=\varepsilon_{p_i}(f_*)\in{\mathbb Z}_+$, $i=1,\dots,r$.
Obviously, the tuple of integers $(b_1,\dots,b_r)\in{\mathbb
Z}^r_+$ does not depend on the choice of the tuple $(f_*)$ and
makes an invariant of the subvariety $B$. We will assume that the
inetegers $b_i$ are ordered: $b_1\geq b_2\geq\dots\geq
b_r$.\vspace{0.1cm}

{\bf Definition 5.2.}  The ordered (non-increasing) tuple of
integers $(b_1,\dots, b_r)$ is called the {\it global type} of the
bi-invariant subvariety $B\subset Y_r$, the notation:
$\varepsilon_{\rm total}(B)=(b_*)$.\vspace{0.1cm}

Set
$$
r_*=\mathop{\rm max}\{j\,|\,b_j\geq j, 1\leq j\leq r\},
$$
if $b_1\geq 1$; if $b_1=\dots=b_r=0$, then we set $r_*=0$. Set
$\Phi(0,\dots,0)=0$, and for $b_1\geq 1$
$$
\Phi(b_1,\dots,b_r)=\sum^{r_*}_{j=1}(b_j+1)(b_j+1-j).
$$
Now we have\vspace{0.1cm}

{\bf Lemma 5.2.} {\it The following estimate holds:}
$$
\mathop{\rm codim}(B\subset Y_r)\geq\Phi(b_1,\dots,b_r).
$$
\vspace{0.1cm}

{\bf Proof.} If $b_1=\dots=b_r=0$, then there is nothing to prove.
So assume that $b_1\geq 1$. By the bi-invariance, it is sufficient
to show that for fixed linearly independent points
$p_1,\dots,p_r\in{\mathbb P}^N$ the estimate
$$
\mathop{\rm codim}(B(p_1,\dots,p_r)\subset
Y_r(p_1,\dots,p_r))\geq\Phi(b_*)
$$
holds, where $Y_r(p_*)=\{(f_*)\in
Y_r\,|\,\{p_1,\dots,p_r\}=|Z(f_*)|\}$ и $B(p_*)=B\cap Y_r(p_*)$.
Furthermore, taking an affine chart ${\mathbb
C}^N_{(z_1,\dots,z_N)}\subset{\mathbb P}^N$, we may assume that
$p_1$ is the origin, and for $j\geq 2$
$$
p_j=(0,\dots,0,1,0,\dots,0),
$$
where the only unity occupies the $(j-1)$-th position. The closed
subset $B(p_*)\subset Y_r(p_*)$ may consist of several irreducible
components. Take any component $B^+$ of that set, the generic
tuple $(f_*)\in B^+$ in which satisfies the equalities
$$
\varepsilon_{p_i}(f_*)=b_i.
$$
Obviously, the lemma will be shown if we prove the inequality
$$
\mathop{\rm codim}(B^+\subset Y_r(p_*))\geq\Phi(b_*).
$$
That is what we will do. To simplify the arguments, assume that
all polynomials $f_i$ are quadratic: if $\mathop{\rm deg}f_i=3$,
then the arguments work without modification, whereas the
estimates for the codimension only get stronger (there are more
coefficients).\vspace{0.1cm}

Write down explicitly
$$
f_i=a^{(i)}_1z_1+\dots+a^{(i)}_Nz_N+\sum_{j\leq k}a^{(i)}_{jk}z_jz_k.
$$
If $r\geq 2$, then the condition $f_i(p_j)=0$ takes the form of
the equalities
$$
a^{(i)}_j+a^{(i)}_{jj}=0
$$
for all $i$ and $j=1,\dots,r-1$. By the bi-invariance, we may
assume that
$$
\mathop{\rm rk}(df_1(o),\dots,df_{N-b_1}(o))=N-b_1,
$$
and the linear forms $df_i(o)$ are linear combinations of the
first $N-b_1$ forms $df_1(o),\dots,df_{N-b_1}(o)$ for $i\geq
N-b_1+1$. This gives $b_1(b_1+1)$ independent conditions on the
coefficients $a^{(i)}_j$ for $i\geq N-b_1+1$, assuming the
polynomials $f_1,\dots,f_{N-b}$ to be fixed. If $r=1$, then there
is nothing more to prove.\vspace{0.1cm}

Assume that $r\geq 2$ and consider the conditions, associated with
the point $p_2=(1,0,\dots,0)$. Recall that $b_2\leq b_1$. If
$r_*=1$, there is nothing to prove. Therefore we assume that
$b_2\geq 2$. Again we assume the polynomials $f_1,\dots,f_{N-b_2}$
to be fixed (which does not contradict the first step of the proof
above), so that the linear forms $df_1(p_2),\dots,df_{N-b_2}(p_2)$
are linearly independent, and $df_i(p_2)$ for $i\geq N-b_2+1$ are
their linear combinations. Explicitly,
$$
df_i(p_2)= (a^{(i)}_1+2a^{(i)}_{11})z_1+\sum_{j\geq
2}(a^{(i)}_j+a^{(i)}_{1j})z_j.
$$
By the equality $a^{(i)}_1=-a^{(i)}_{11}$ the coefficient at $z_1$
is linearly dependent on the set of coefficients $a^{(i)}_j$.
However, the coefficients $a^{(i)}_{1j}$ were not involved in the
conditions, associated with the point $p_1$. Therefore, requiring
that
\begin{equation}\label{25april12.1}
df_i(p_2)|_{\{z_1=0\}}\in<df_1(p_2)|_{\{z_1=0\}},\dots,df_{N-b_2}(p_2)|_{\{z_1=0\}}>
\end{equation}
for $i\geq N-b_2+1$, we impose at least
$$
(b_2+1)(b_2-1)
$$
independent conditions on the coefficients $a^{(i)}_{1j}$. (The
precise number of conditions is determined by the dimension of the
space in the right hand side of the formula (\ref{25april12.1}):
if it is equal to $N-b_2$, then we get $(b_2+1)(b_2-1)$
conditions, if it drops by one, then we get $(b_2+1)b_2$
conditions.) This completes our consideration of the second
component $(j=2)$ of the function $\Phi$. If $r_*=2$, our lemma is
shown.\vspace{0.1cm}

If $r_*\geq 3$, then we continue in the same spirit: arguing by
induction, we assume that it is shown that the condition
$$
\mathop{\rm rk}(df_1(p_{\alpha}),\dots,(df_{N+1}(p_{\alpha}))=N-b_{\alpha}
$$
for $\alpha=1,\dots,j$ imposes on the coefficients
$a^{(i)}_\gamma$, $i\geq N-b_1+1$, and $a^{(i)}_{\lambda k}$,
where $\lambda=1,\dots,j-1$, $k=\lambda+1,\dots,N$ and $i\geq
N-b_{\lambda+1}+1$, in total at least
$$
\sum^j_{\lambda=1}(b_{\lambda}+1)(b_{\lambda}+1-\lambda)
$$
independent conditions. If $r_*=j$, then the proof is complete at
that step. Otherwise, assuming $f_1,\dots,f_{N-b_{j+1}}$ to be
fixed, we get that the linear forms
$$
df_i(p_{j+1})|_{\{z_1=\dots=z_j=0\}},\quad i\geq N-b_{j+1}+1,
$$
belong to the linear space
$$
<df_1(p_{j+1})|_{\{z_1=\dots=z_j=0\}},\dots,
df_{N-b_{j+1}}(p_{j+1})|_{\{z_1=\dots=z_j=0\}}>.
$$
This gives at least $(b_{j+1}+1)(b_j-j)$ independent conditions
for the {\it new} (that is, not involved in the previous
considerations) coefficients
$$
a^{(i)}_{j,j+1},\dots,a^{(i)}_{j,N},\quad i\geq N-b_{j+1}+1.
$$
Now the inductive step from $j$ to $j+1$ is constructed and the
proof of the lemma is complete.\vspace{0.1cm}

{\bf Corollary 5.1.}  {\it The following estimate holds:}
$$
\mu_{\rm total}(i,a,d)\leq \sum_{\Phi(b_1,\dots,b_i)+i\leq
a}\mu_{\rm local}(a-b_i,b_i;d).
$$
\vspace{0.1cm}

{\bf Proof.} This follows from the previous lemma, the equality
$\mathop{\rm codim}Y_i=i$ and the fact that the global
multiplicity $\mu_{\rm total}$ is computed via the tuples of $N+1$
polynomials, whereas the local one $\mu_{\rm local}$ via the
tuples of $N$ polynomials. Indeed, if
$$
B\subset{\cal P}^{\times n_1}_{[1,2],N}\times{\cal P}^{\times n_2}_{[1,3],N}
$$
is a bi-invariant subvariety of codimension $a$ (in the sense of
the local theory \S\S2-4) and $\varepsilon(B)=b$, where
$n_1+n_2=N+1$, then the projection $[B]_{N+1}$ of the set $B$
along the last direct factor has codimension at most $a-b$ in the
space ${\cal P}^{\times n_1}_{[1,2],N}\times{\cal P}^{\times
(n_2-1)}_{[1,3],N}$, since for a generic tuple
$(f_1,\dots,f_{N+1})\in B$ we have: $df_{N+1}(o)$ vanishes on a
subspace of codimension $b$, that depends on $(f_1,\dots,f_N)$
only. Now the claim of the corollary is obvious.
Q.E.D.\vspace{0.1cm}

%%%%%%%%%%%%%%%%%%%%%%%%%%%%%%%%%%%%%%%%%%%%%%%%%%%%%%%%%%%%%%%%%
%%%%%%%%%%%%%%%%%%%%%%%%%%%%%%%%%%%%%%%%%%%%%%%%%%%%%%%%%%%%%%%%%
%%%%%%%%%%%%%%%%%%%%%%%%%%%%%%%%%%% subsection 5.4

To simplify the notations, we will omit the parameters $l$ and
$N$, the more so that our estimates do not depend on those
parameters, and write
$$
\mu_{\rm local}(a, b; d)
$$
instead of $\mu^{\rm local}_{l,N}(a,b;d)$. The values of the
parameters $l$ and $N$ are in any case fixed in the subsequent
arguments.\vspace{0.3cm}

{\bf 5.4. An explicit estimate for the global multiplicity.} Now
everything is ready to obtain an effective estimate for the global
multiplicity $\mu_{\rm total}(r,a,d)$. Let $B\subset Y_r$ be an
irreducible bi-invariant subvariety, the codimension of which in
the space ${\cal H}^{(n_1,n_2)}_{N+1}$ does not exceed $a$. In
particular, the inequality $\Phi(b_1,\dots,b_r)+r\leq a$ holds,
where $(b_1,\dots,b_r)=\varepsilon(B)$ is the global type of the
subvariety $B$.\vspace{0.1cm}

{\bf Proposition 5.4.} {\it The following inequality holds:}
$$
\mu_{\rm total}(B,d)\leq \sum^r_{i=1}\mu_{\rm local}(a-b_i,b_i;d).
$$
\vspace{0.1cm}

{\bf Corollary 5.2.}  (i) {\it The following inequality holds:}
$$
\mu_{\rm total}(r,a,d)\leq\mathop{\rm max}_{\Phi(b_1,\dots,b_r)+r\leq a}
\left\{\sum^r_{i=1}\mu_{\rm local}(a-b_i,b_i;d)\right\}.
$$
\vspace{0.1cm}

(ii) {\it The following inequality holds:}
$$
\mu_{\rm total}(r,a,d)\leq d\cdot\left(\mathop{\rm max}_{\Phi(b_1,\dots,b_r)+r\leq a}
\left\{\sum^r_{i=1}\bar{\mu}_{\rm local}(a-b_i,b_i)\right\}\right).
$$
\vspace{0.1cm}

{\bf Proof of the corollary.} The first inequality follows
immediately from Proposition 5.4 by the definition of the numbers
$\mu_{\rm total}(r,a,d)$. The claim (ii) follows from (i), taking
into account Proposition 4.1. Q.E.D. for the
corollary.\vspace{0.1cm}

{\bf Proof of Proposition 5.4.} Let $(f_*)\in B$ be a tuple of
general position, $\{p_1,\dots,p_r\}=|Z(f_*)|$ its common zeros,
where $b_i=\varepsilon_{p_i}(f_*)$, $i=1,\dots,r$. The inequality
of Proposition 5.4 follows from the estimate
\begin{equation}\label{25april12.2}
\mu_{\rm local}((f_*); R,p_i)\leq\mu_{\rm local}(a-b_i,b_i;d)
\end{equation}
for every effective cycle $R$ of pure codimension $l$ and degree
$d$. It is the last estimate that we will prove.\vspace{0.1cm}

Set $p=p_i$ and let $B(p)\subset B$ be the closed (in $B$) subset
of tuples $(g_*)$, vanishing at the point $p$. By the
bi-invariance of the set $B$ the original tuple $(f_*)$, fixed at
the beginning of the proof, is a generic tuple of one of the
irreducible components $B^+$ of the set $B(p)$, and moreover, the
codimension of $B^+$ in the space ${\cal P}^{(n_1,n_2)}_N$
coincides with $\mathop{\rm codim}(B\subset{\cal
H}^{(n_1,n_2)}_{N+1})$ and for that reason does not exceed $a$ (we
used the natural identification of the linear space of tuples
$(g_*)\in{\cal H}^{(n_1,n_2)}_{N+1}$, vanishing at the point $p$,
with the space ${\cal P}^{(n_1,n_2)}_N$, defined in \S2). Let
$$
\pi\colon{\cal P}^{(n_1,n_2)}_N\to{\cal P}^{(n_1,n_2-1)}_N
$$
be the projection along the last direct factor and $[B^+]_N$ the
closure of the image $\pi(B^+)$. By the bi-invariance and the
condition $\varepsilon_p(f_*)=b$, the generic fibre of the
projection
$$
\pi_B\colon B^+\to[B^+]_N
$$
has in ${\cal P}_{[1,3],N}$ codimension at least $b$
($df_{N+1}(p)$ vanishes on a $b$-dimensional linear subspace that
depends on $df_1,\dots,df_N$ only). Therefore, the inequality
$$
\mathop{\rm codim}([B^+]_N\subset{\cal P}^{(n_1,n_2-1)}_N)\leq a-b
$$
holds, so that the inequality (\ref{25april12.2}) is shown. Q.E.D.
for Proposition 5.4.\vspace{0.1cm}

Set
$$
\mu_{\rm total}(a,d)=\mathop{\rm max}_{1\leq r\leq {\rm min}\{a,N+1\}}
\{\mu_{\rm total}(r,a,d)\}.
$$
\vspace{0.1cm}

{\bf Proposition 5.5.} {\it For $a\geq 12$ the following estimate
holds:}
$$
\mu_{\rm total}(a,d)\leq 3\cdot 2^{a-6}d
$$
\vspace{0.1cm}

{\bf Proof.} Assume first that $a\geq 21$. Since $r\leq a$, from
the claim (i) of Corollary 5.2 and the inequality (i) of Corollary
4.3 we get
$$
\mu_{\rm total}(a,d)\leq
a\frac{e^2}{2\pi[\sqrt{a}]}\left(\frac53\right)^{[\sqrt{a}]}d.
$$
It is easy to check that for $a\geq 21$ the right hand side of the
last inequality is strictly smaller than $3\cdot 2^{a-6}$. This
method, however, is very coarse. A mush more precise estimate is
given by the claim (ii) of corollary 5.2. Set
$$
\bar{\mu}_{\rm total}(a)=\mathop{\rm max}_{1\leq r\leq a}\left(
\mathop{\rm max}_{\Phi(b_*)+r\leq a}\left\{
\sum^r_{i=1}\bar{\mu}(a-b_i,b_i)\right\}\right).
$$
According to Corollary 5.2, (ii), and Proposition 4.1 the
inequality
$$
\mu_{\rm total}(a,d)\leq d\bar{\mu}_{\rm total}(a)
$$
holds, so that to complete the proof of Proposition 5.5, it is
sufficient to check that for $12\leq a\leq 20$ the inequality
$\bar{\mu}_{\rm total}(a)\leq 3\cdot 2^{a-6}$ holds. For small
values of the codimension $a$ the function $\bar{\mu}_{\rm
total}(a)$ is easy to compute explicitly, and the results (for
$1\leq a\leq 36$) are given below in the table, which is organized
in the following way. Each row corresponds to a certain value of
$a$, which is given in the first column. In the second column of
the same row we give the values of the parameters $r$ (the number
of points) and $b_1,\dots,b_r$, for which the maximum in the
definition of the function $\bar{\mu}_{\rm total}$ is attained. In
the third column we give the very value of $\bar{\mu}_{\rm
total}(a)$.\vspace{0.5cm}

\begin{tabular}{|l|l|l|}
\hline
$a$   &   &   $\bar{\mu}_{\rm total}(a)$ \\
\hline
1     &   $r=1$, $b_1=0$          & 1     \\
\hline
2     &   $r=2$, $b_1=b_2=0$          & 2     \\
\hline
3     &   $r=3$, $b_1=b_2=b_3=0$          & 3     \\
\hline
4     &   $r=4$, $b_1=b_2=1$          & 6     \\
\hline
5     &   $r=3$, $b_1=b_2=b_3=1$          & 15     \\
\hline
6     &   $r=4$, $b_1=\dots=b_4=1$          & 24     \\
\hline
7     &   $r=5$, $b_1=\dots=b_5=1$          & 35     \\
\hline
8     &   $r=6$, $b_1=\dots=b_6=1$          & 48     \\
\hline
9     &   $r=7$, $b_1=\dots=b_7=1$          & 63     \\
\hline
10     &   $r=8$, $b_1=\dots=b_8=1$          & 80     \\
\hline
11     &   $r=9$, $b_1=\dots=b_9=1$          & 99     \\
\hline
12     &   $r=10$, $b_1=\dots=b_{10}=1$          & 120     \\
\hline
13     &   $r=11$, $b_1=\dots=b_{11}=1$          & 143     \\
\hline
14     &   $r=12$, $b_1=\dots=b_{12}=1$          & 168     \\
\hline
15     &   $r=13$, $b_1=\dots=b_{13}=1$          & 195     \\
\hline
16     &   $r=7$, $b_1=\dots=b_{7}=2$          & 308     \\
\hline
17     &   $r=8$, $b_1=\dots=b_{8}=2$          & 408     \\
\hline
18     &   $r=9$, $b_1=\dots=b_{9}=2$          & 522     \\
\hline
19     &   $r=10$, $b_1=\dots=b_{10}=2$          & 660     \\
\hline
20     &   $r=11$, $b_1=\dots=b_{11}=2$          & 814     \\
\hline
21     &   $r=12$, $b_1=\dots=b_{12}=2$          & 996     \\
\hline
22     &   $r=13$, $b_1=\dots=b_{13}=2$          & 1196     \\
\hline
23     &   $r=14$, $b_1=\dots=b_{14}=2$          & 1428     \\
\hline
24     &   $r=15$, $b_1=\dots=b_{15}=2$          & 1680     \\
\hline
25     &   $r=16$, $b_1=\dots=b_{16}=2$          & 1968     \\
\hline
26     &   $r=17$, $b_1=\dots=b_{17}=2$          & 2278     \\
\hline
27     &   $r=18$, $b_1=\dots=b_{18}=2$          & 2628     \\
\hline
\end{tabular}
\vspace{0.5cm}

The subsequent values of the function $\bar{\mu}_{\rm total}$ are
as follows:
$$
\bar{\mu}_{\rm total}(28)=3002,\quad\bar{\mu}_{\rm total}(29)=3420,
\quad\bar{\mu}_{\rm total}(30)=3864,
$$
$$
\bar{\mu}_{\rm total}(31)=4356,\quad\bar{\mu}_{\rm total}(32)=4876,
\quad\bar{\mu}_{\rm total}(33)=5448,
$$
$$
\bar{\mu}_{\rm total}(34)=6050,\quad\bar{\mu}_{\rm total}(35)=6708.
$$
All these values are attained for $r=a-9$, $b_1=\dots=b_r=2$. For
$a=36$ the next jump takes place: the maximum is attained for
$r=12$ and $b_1=\dots=b_{12}=3$ and equal to $\bar{\mu}_{\rm
total}(36)=7980$.\vspace{0.1cm}

According to the table given above, starting from $a=12$ the
required inequality $\bar{\mu}_{\rm total}(a)\leq 3\cdot 2^{a-6}$
holds.\vspace{0.1cm}

This proves the proposition. Q.E.D.\vspace{0.1cm}

Finally, from the values of the function $\bar{\mu}_{\rm total}$,
given in the table, by elementary arithmetic we get\vspace{0.1cm}

{\bf Proposition 5.6.} {\it If the triple $(a,n_1,n_2)$ is one of
the following triples:
$$
(11,5,3),\quad(11,3,4),\quad(11,1,5),\quad (10,2,4),\quad(10,0,5),
$$
then the following inequality holds:}
$$
\bar{\mu}_{\rm total}(a,d)\leq 2^{n_1+n_2-4}3^{n_2-1}d.
$$
\vspace{0.3cm}

%%%%%%%%%%%%%%%%%%%%%%%%%%%%%%%%%%%%%%%%%%%%%%%%%%%%%%%%%%%%%%%%%%%%%%%%%%%%%%%%%%%%
%%%%%%%%%%%%%%%%%%%%%%%%%%%%%%%%%%%%%%%%%%%%%%%%%%%%%%%%%%%%%%%%%%%%%%%%%%%%%%%%%%%%
%%%%%%%%%%%%%%%%%%%%%%%%%%%%%%%%%%% subsection 5.5

{\bf 5.5. Regular complete intersections.} Finally, let us prove
Theorem 2. Let
$$
{\cal B}\subset{\mathbb P}\times\left({\cal P}^{\times k_1}_{2,M+k+1}
\times {\cal P}^{\times k_2}_{3,M+k+1}\right)
$$
be the closed set of ``bad'' pairs $(o,(f_*))$, where
$f_1(o)=\dots=f_k(o)=0$ and at least one of the conditions (R1-R3)
is violated at this point. Let $\pi_1$ and $\pi_2$ be the
projections of the direct product on the first $({\mathbb P})$ and
second (the space of tuples of $k$ polynomials) factors,
respectively. Set ${\cal B}(o)=\pi^{-1}_1(o)\cap{\cal B}$. The
closed subset ${\cal B}(o)$ is contained in the subspace
$$
{\cal L}(o)={\cal P}^{\times k_1}_{[1,2],M+k}\times {\cal
P}^{\times k_2}_{[1,3],M+k} \subset\pi^{-1}_1(o)
$$
(where we again identify homogeneous polynomials vanishing at the
point $o$, with non-homogeneous polynomials without the free
term), and it is sufficient to show that its codimension with
respect to that subspace is at least $M+1$. Indeed, if this is the
case, then
$$
\mathop{\rm codim}{\cal B}=\mathop{\rm codim}({\cal B}(o)\subset\pi^{-1}_1(o))\geq M+k+1,
$$
so that the map $\pi_2|_{\cal B}$ can not be surjective, which
immediately implies the claim of Theorem 2.\vspace{0.1cm}

For the conditions (R1) and (R2) the inequality
\begin{equation}\label{25april12.3}
\mathop{\rm codim}({\cal B}(o)\subset{\cal L}(o))\geq M+1
\end{equation}
is shown in \cite{Pukh01}, taking into account Proposition 5.2.
Thus it is sufficient to prove the inequality (\ref{25april12.3})
for such tuples $(f_*)\in{\cal B}(o)$, which satisfy the
conditions (R1) and (R2), but not the condition (R3). Fixing the
linear parts of the polynomials $f_i$ at the point $o$, and thus
the projectivized tangent space ${\mathbb T}\cong{\mathbb
P}^{M-1}$, we reduce the problem to estimating the codimension of
the set of tuples of $n_1=k_1+k_2=k$ quadratic and $n_2=k_2$ cubic
homogeneous polynomials
\begin{equation}\label{25april12.4}
(\bar{q}_{1,2},\dots,\bar{q}_{k,2},\bar{q}_{k_1+1,3},\dots,\bar{q}_{k,3}),
\end{equation}
for which the condition (R3) is not satisfied. However, by
Propositions 5.5 and 5.6 the main inequality (\ref{9april12.7}) of
the condition (R3) is satisfied for a generic element of every
subvariety of codimension $\leq M$ in the space of tuples
$$
{\cal P}^{\times n_1}_{2,M}\times{\cal P}^{\times n_2}_{3,M}
$$
and every irreducible subvariety $R\subset{\mathbb T}$ of
codimension three. Therefore, the closed set of tuples
(\ref{25april12.4}), not satisfying the condition (R3), is of
codimension at least $M+1$. This proves the estimate
(\ref{25april12.3}) and Theorem 2 as well. Q.E.D.\vspace{0.1cm}

%%%%%%%%%%%%%%%%%%%%%%%%%%%%%%%%%%%%%%%%%%%%%%%%%%%%%%%%%%%%%%%%%%%%%%%%%%%%%%%%%%%%
%%%%%%%%%%%%%%%%%%%%%%%%%%%%%%%%%%%%%%%%%%%%%%%%%%%%%%%%%%%%%%%%%%%%%%%%%%%%%%%%%%%%
%%%%%%%%%%%%%%%%%%%%%%%%%%%%%%%%%%% BIBLIOGRAPHY

\begin{flushleft}
Department of Mathematical Sciences,\\
The University of Liverpool
\end{flushleft}

\noindent{\it pukh@liv.ac.uk}

\end{document}